\documentclass[12pt]{article}

\pdfoutput = 1

\usepackage[english]{babel}
\usepackage[latin1]{inputenc}

\usepackage{geometry}
\usepackage{fancyhdr}
\usepackage{setspace}
\newlength{\separaParrafos}
\parskip = \separaParrafos	
\usepackage{float}

\usepackage[
	bookmarksnumbered,
	colorlinks,
	linkcolor = black
	]{hyperref}

\usepackage{graphicx}
\usepackage{tabularx}
\usepackage{amsfonts}
\usepackage{amssymb}
\usepackage{pstricks}
\usepackage{amsmath,mathtools,upgreek}
\usepackage{array}
\usepackage{pgfplots}
\usepackage{cancel}
\usepackage{bm}

\usepackage{caption}
\usepackage{subcaption}

\pgfplotscreateplotcyclelist{ageplot}{%
{blue,mark=o, mark size={3}},
{red,mark=square},
{green,mark=star},
{brown,mark=triangle},
{black,mark=diamond},}

\pgfplotsset{compat=1.3}
\pgfplotsset{width=0.5\textwidth}
\pgfplotsset{
tick label style={font=\small},
label style={font=\small},
legend style={font=\small},every axis/.append style={line width=1pt,tick
style={semithick} }
}
\pgfplotstableset{col sep=comma,row sep=newline}

\newcommand{\vm}[1]{\bm{\mathrm{#1}}} 		
\newcommand{\norm}[1]{\lVert#1\rVert} 		
\newcommand{\abs}[1]{\left|#1\right|} 		
\newcommand{\sst}{\vm{\upsigma}^{*}}		
\newcommand{\bb}{\vm{b}}			
\newcommand{\bsigma}{\vm{\upsigma}}		
\newcommand{\bv}{\vm{v}}			
\newcommand{\abil}[2]{a(#1,#2)}			
\newcommand{\babil}[2]{\bar{a}(#1,#2)}		
\newcommand{\llin}[1]{l({#1})}			
\newcommand{\V}{V}				
\newcommand{\bu}{\vm{u}}			
\newcommand{\uH}{u^H}				
\newcommand{\be}{\vm{e}}			
\newcommand{\bn}{\vm{n}}			

\title{Report:\\
Error estimation of recovered solution in FE analysis}

\author{E. Nadal$^1$, O.A. González-Estrada$^2$, J.J. Ródenas$^1$, F.J. Fuenmayor$^1$}

\date{}

\begin{document}

\maketitle

\begin{center}\small{
$^{1}$Centro de Investigaci\'on de Tecnolog\'ia de Veh\'iculos (CITV), \\Universitat Polit\`{e}cnica de Val\`{e}ncia, E-46022-Valencia, Spain, \\e-mail: ennaso@upvnet.upv.es, jjrodena@mcm.upv.es, ffuenmay@mcm.upv.es\\
$^{2}$Institute of Mechanics and Advanced Materials (IMAM), Cardiff School of Engineering, Cardiff University, 
Queen's Buildings, The Parade, Cardiff CF24 3AA Wales, UK, e-mail: estradaoag@cardiff.ac.uk}
\end{center}

\textbf{Abstract.} The recovery type error estimators introduced by Zienkiewicz and Zhu use a recovered stress field evaluated from the Finite Element (FE) solution. Their accuracy depends on the quality of the recovered field. In this sense, accurate results are obtained using recovery procedures based on the Superconvergent Patch Recovery technique (SPR). These error estimators can be easily implemented and provide accurate estimates. Another important feature is that the recovered solution is of a better quality than the FE solution and can therefore be used as an enhanced solution.

We have developed an SPR-type recovery technique that considers equilibrium and displacements constraints to obtain a very accurate recovered displacements field from which a recovered stress field can also be evaluated. We propose the use of these recovered fields as the standard output of the FE code instead of the raw FE solution. Techniques to quantify the error of the recovered solution are therefore needed. 

In this report we present an error estimation technique that accurately evaluates the error of the recovered solution both at global and local levels in the FEM and XFEM frameworks. We have also developed an \textit{h}-adaptive mesh refinement strategy based on the error of the recovered solution. As the converge rate of the error of the recovered solution is higher than that of the FE one, the computational cost required to obtain a solution with a prescribed accuracy is smaller than for traditional \textit{h}-adaptive processes.


%
%
\section{Introduction}
%
%
Here we briefly discuss the possibility of estimating the error in the recovered
solution $\vm{\sigma}^*$ evaluated from the FE solution $\vm{\sigma}^h$ as part of the process of error estimation in energy norm based on  recovery-techniques. Recovery error estimation techniques are based on the assumption that the recovered solution is more accurate than the FE solution. A sufficiently accurate estimation of the error in energy norm of the recovered solution could lead to refinement processes based on the accuracy of the recovered solution instead of that of the FE solution. This could result in a considerable reduction of the computational cost that would be particularly interesting, for example, in optimization processes whose efficiency would be significantly increased as the required level of accuracy would be gained with a considerable lower number of degrees of freedom.

In \cite{diezrodenas2007} and the Ph.D. thesis \cite{gonzalez2010} we developed techniques whose objective was that of obtaining upper bounds of the error in energy norm by means of the use of recovery-based error estimators. This 
research yielded an expression of an error bounding technique which includes two types of terms. The first term is the expression of the ZZ error estimator. The second term is a correction term used to account for the lack of equilibrium of the recovered stress field (both in the interior of the domain and along the domain boundaries). Numerical experiments have shown that slightly manipulating the expression of the correction term one can estimate the error of the recovered solution, not only at a global level but also at element level. 

In this report we will show the initial developments that finally led to an empirical expression based on this correction term that can be used to efficiently estimate the error in energy norm of the recovered solution $\vm{u}^*$, initially evaluated to estimate the error of the FE solution $\vm{u}^h$. So, the report is structured  following our developments along time. Section \ref{sec:ModProb} shows the model problem. Section \ref{sec:UpperBounds} shows the technique used to evaluate upper bounds of the error based on nearly equilibrated recovered solutions that leads to the correction term previously mentioned. Section \ref{sec:URecovery} describes the technique to obtain the recovered displacement and stress fields ($\vm{u}^*$ and $\vm{\sigma}^*$) that will be used to estimate the discretization error in energy norm of the finite element solution $\vm{u}^h$. Section \ref{sec:ErrRecovered} will present an expression that could be used to evaluate the error in energy norm of the recovered solution and the empirical approximation to this expression that can be used in practice to efficiently estimate this error. Finally, the numerical examples showing the potential of the technique proposed for the estimation of the error in the recovered solution are presented in Section \ref{sec:Examples}.
%
%
\section{Model problem} 
\label{sec:ModProb}
%
%
Let us consider the 2D linear elastic problem. The unknown displacement field
$\vm{u}$, taking values in $\Omega \subset \mathbb{R}^{2}$, is the solution of
the boundary value problem given by 
\begin{align}
  -\nabla \cdot \vm{\sigma} \left(\vm{u}\right) &= \vm{b}  	&&  {\rm in }\;
\Omega 
   \label{Eq:IntEq} \\
   \vm{\sigma} \left(\vm{u} \right)\cdot \vm{n} &= \vm{t} 	&&  {\rm on }\;
\Gamma _{N}  		\label{Eq:Neumann}\\
   \vm{u}                                         &= \vm{0}	&&  {\rm on }\;
\Gamma _{D} \label{Eq:Dirichlet}
\end{align}

\noindent where $\Gamma _{N}$ and $\Gamma _{D}$ denote the Neumann and Dirichlet
boundaries with $\partial \Omega = \Gamma_N \cup \Gamma_D$ and $\Gamma_N \cap
\Gamma_D =\varnothing$. The Dirichlet boundary condition in (\ref{Eq:Dirichlet})
is taken homogeneous for the sake of simplicity. 
 
The weak form of the problem reads: Find $\vm{u} \in V$ such that 
\begin{equation} \label{Eq:WeakForm} 
a(\vm{u},\vm{v}) = l(\vm{v}) \qquad \forall \vm{v} \in V, 
\end{equation}  
 
\noindent where ${V}$ is the standard test space for the elasticity problem
such that $V = \{\vm{v} \;|\; \vm{v} \in  H^1(\Omega) ,
\vm{v}|_{\Gamma_D}(\vm{x}) = \vm{0} \}$, and 
\begin{align}
a(\vm{u},\vm{v})  & := \int _{\Omega}  \vm{\sigma} (\vm{u}): 
\vm{\varepsilon}(\vm{v}) d \Omega = 
\int _{\Omega}  \vm{\sigma}(\vm{u}) : \vm{\mathsf{S}} : \vm{\sigma} (\vm{v}) d
\Omega \\
l(\vm{v})&:=\int _{\Omega} \vm{b}  \cdot  \vm{v}d \Omega + \int _{\Gamma_N}
\vm{t}  \cdot \vm{v}d \Gamma,
\end{align}
 
\noindent where $ \vm{\sigma}$ and $\vm{\varepsilon}$ denote the stresses and
strains, and  $\vm{\mathsf{S}}$ is the compliance tensor.

The bilinear form $a(\cdot,\cdot)$ can also be expressed in terms of stresses by formally introducing $\bar{a}(\cdot,\cdot)$ such that $\bar{a}(\vm{\sigma}, \vm{\tau}):=\int_{\Omega}{\vm{\sigma}:\vm{D}^{-1}\vm{\tau}d\Omega}$. Note that $a(\vm{u},\vm{v})=\bar{a}(\vm{\sigma},\vm{\tau})$.

Let $\vm{u}^{h}$ be a finite element approximation to  $\vm{u}$. The solution
lies in a functional space $V^{h} \subset V$ associated with a mesh of finite elements of characteristic size $h$, and it is such that 
\begin{equation}  
a(\vm{u}^{h},\vm{v}) = l(\vm{v}) \qquad  \forall \vm{v} \in V^{h} 
\end{equation}  

We focus on assessing the error $\vm{e} =\vm{u} -\vm{u}^h$. The error is
measured in the energy norm, induced by $a(\cdot,\cdot)$. The quantity to be
assessed is $\norm{\vm{e}}^2 = a(\vm{e},\vm{e})$.
%
%
\section{Upper error bound  of the FE solution from nearly equilibrated recovered solutions}
\label{sec:UpperBounds}
%
%
This section will explain with detail the basis of the proposed upper error bounding method. The ideas were presented in \cite{diezrodenas2007} and the Ph.D. thesis \cite{gonzalez2010}. 

\noindent\textbf{Proposition.} Let $\sst$ be a recovered field such that $\sst \cdot
\vm{n}$ is continuous almost everywhere along any interior curve
$\Gamma \subset \Omega$ (being $\vm{n}$ the outward normal vector to $\Gamma$)
and
\begin{subequations}\label{eq:strong}
\begin{align}
 - \nabla \cdot \sst & = \bb + \vm{s} & \text{ a.e. in } \Omega \label{eq:EDP}
\\[2ex]
 \sst \cdot \vm{n} & = \vm{t} +  \vm{r} & \text{ a.e. on } \Gamma_N
\label{eq:Neumann} 
\end{align}
\end{subequations}

Note that if $\vm{s}=0$ and $\vm{r}=0$, $\sst$ is statically
admissible. For \emph{small} $\vm{s}$ and $\vm{r}$ (when compared to
$\vm{b}$ and $\vm{t}$, respectively) $\sst$ is said to be nearly
statically admissible, up to the equilibrium default terms.

Then, the following expression holds:
\begin{equation}\label{eq:w1}
\babil{\sst}{\bsigma(\bv)} = \llin{\bv}  + \int_\Omega
\bv\cdot\vm{s}\,d\Omega + \int_{\Gamma_N} \bv\cdot\vm{r}\,d\Gamma
\;\; \forall \bv \in \V,
\end{equation}
and, as a direct consequence assuming that the Dirichlet boundary conditions are
homogeneous
\begin{equation}\label{eq:w2}
\norm{\bu}^2=\babil{\bsigma(\bu)}{\bsigma(\bu)}\le
\babil{\sst}{\sst}-2\int_\Omega
\bu\cdot\vm{s}\,d\Omega-2\int_{\Gamma_N} \bu\cdot\vm{r}\,d\Gamma
\end{equation}

Moreover, the error approximation $\sst_e=\sst-\bsigma(\vm{u}^h)$ could also lead to an upper bound of the error norm:  

\begin{equation}\label{eq:w3}
\norm{\be}^2=\abil{\be}{\be}=\babil{\bsigma_e}{\bsigma_e}\le
\babil{\sst_e}{\sst_e}-2\int_\Omega
\vm{e}\cdot\vm{s}\,d\Omega-2\int_{\Gamma_N}
\vm{e}\cdot\vm{r}\,d\Gamma
\end{equation}
and this latter property stands also for non-homogeneous
Dirichlet boundary conditions, up to oscillation terms.

\noindent\textbf{Proof.}
The expression in \eqref{eq:w1} is derived from \eqref{eq:strong} using
the standard weighted residuals technique combined with
integration by parts  (and the fact that the test function $\bv$
vanishes on $\Gamma_N$)
\begin{align*}
\int_\Omega \bv\cdot(\bb + \vm{s}) \, d\Omega &=
-\int_\Omega \bv \cdot \nabla \cdot \sst \,d\Omega\\
&=- \int_\Omega  \nabla \cdot ( \sst \cdot  \bv) \,d\Omega + \int_\Omega  \nabla
\bv :  \sst \,d\Omega\\
&=-\int_{\partial \Omega}   \bn \cdot (\sst \cdot  \bv) \,d\Gamma+ \int_\Omega 
\frac{1}{2} (\nabla \bv +  \nabla^{\textsf{T}} \bv )\cdot  \sst \,d\Omega\\
&=-\int_{\Gamma_N}   \bv \cdot (\sst \cdot  \bn) \,d\Gamma+
\babil{\sst}{\bsigma(\bv)}\\
&=-\int_{\Gamma_N}   \bv \cdot (\vm{t} +  \vm{r}) \,d\Gamma+
\babil{\sst}{\bsigma(\bv)}
\end{align*}

Then \eqref{eq:w1} follows using the definition of $\llin{\bv}$. Thus, the proof of \eqref{eq:w2} is straightforward considering in \eqref{eq:w1} $\bv=\bu$
(this requires that in the original problem, the Dirichlet boundary conditions states that $\bu =  \vm{0}$ on $\Gamma_D$) 
\begin{align}\label{eq:p1}
\babil{\sst}{\bsigma(\bu)} &= \llin{\bu}  + \int_\Omega
\bu\cdot\vm{s}\,d\Omega + \int_{\Gamma_N} \bu\cdot\vm{r}\,d\Gamma \nonumber\\
&=\babil{\bsigma(\bu)}{\bsigma(\bu)} + \int_\Omega
\bu\cdot\vm{s}\,d\Omega + \int_{\Gamma_N} \bu\cdot\vm{r}\,d\Gamma
\end{align}
and therefore
\begin{align*}
0\le \babil{\bsigma(\bu)-\sst}{\bsigma(\bu)-\sst}
=&\babil{\bsigma(\bu)}{\bsigma(\bu)}-2\babil{\sst}{\bsigma(\bu)}+\babil{\sst}{
\sst}\\
=&-\babil{\bsigma(\bu)}{\bsigma(\bu)}-2\int_\Omega 
\bu\cdot\vm{s}\,d\Omega \\ 
&-2 \int_{\Gamma_N} \bu\cdot\vm{r}\,d\Gamma+\babil{\sst}{\sst}.\\
\end{align*}

The proof of \eqref{eq:w3} is similar, by considering in \eqref{eq:w1} $\bv=\be$
(this can be done, up to oscillation terms, even if in the original problem the
Dirichlet boundary
conditions are non-homogeneous). Thus,
\begin{align*}\label{eq:p2}
\babil{\sst}{\bsigma(\be)} = \llin{\be}  + \int_\Omega
\be\cdot\vm{s}\,d\Omega + \int_{\Gamma_N} \be\cdot\vm{r}\,d\Gamma
=\babil{\bsigma(\be)}{\bsigma(\be)} \\ 
+\babil{\bsigma(\uH)}{\bsigma(\be)}  + \int_\Omega
\be\cdot\vm{s}\,d\Omega + \int_{\Gamma_N} \be\cdot\vm{r}\,d\Gamma
\end{align*}
that results in
\begin{equation*}\label{eq:p3}
\babil{\sst_e}{\bsigma(\be)} =
\babil{\bsigma(\be)}{\bsigma(\be)} + \int_\Omega
\be\cdot\vm{s}\,d\Omega + \int_{\Gamma_N} \be\cdot\vm{r}\,d\Gamma .
\end{equation*}
Note that the following error representation has been used
\[\babil{\bsigma(\be)}{\bsigma(\be)}=\llin{\be}-\babil{\bsigma(\uH)}{
\bsigma(\be)}.\]

Thus, \eqref{eq:w3} is proved using the same idea as before by simply
considering the positiveness
of $\babil{\bsigma(\be)-\sst_e}{\bsigma(\be)-\sst_e}$.
%
\section{Displacements and stress recovery}
\label{sec:URecovery}
%
%
The evaluation of the correction terms in (\ref{eq:w3}) involves the evaluation of the equilibrium defaults $\vm{s}$ and $\vm{t}$ and  the exact error in the displacements field $\vm{e}$. The evaluation of $\vm{s}$ and $\vm{t}$ depends on the recovery technique used to obtain $\vm{\sigma}^*$ and we are able to derive expressions for the exact evaluations of these terms. However, $\vm{e}=\vm{u}-\vm{u}^h$ can only be estimated. We estimate $\vm{e}$ as $\vm{e} \approx \vm{u}^*-\vm{u}^h$.

We have derived two alternatives for the evaluation of $\vm{e}$. The first one consists in using a sequence of $N$ meshes for the problem. The solution of the last mesh of the sequence would be considered as an enhanced solution $\vm{u}^*$ such that in mesh $i$  
\begin{equation}
	\vm{e}_{(i)}=\vm{u}-\vm{u}^h_{(i)} \approx \vm{u}^h_{(N)}-\vm{u}^h_{(i)}, \ \ i=1,...,N
\label{eq:eAprox}
\end{equation}

In \cite{diezrodenas2007} we evaluated the correction term for the lack of internal equilibrium considering 
\begin{equation}
ES:=-2\int_\Omega \vm{e}\cdot\vm{s} d\Omega\leq 
 2 \left| \int_\Omega \vm{e}\cdot\vm{s} d\Omega \right|  \leq
2\left|\vm{e}\right|_{L_2} \left|\vm{s}\right|_{L_2} 
\label{eq:esCS}
\end{equation}

As previously mentioned, we are able to evaluate $\vm{s}_{(i)}$ and therefore $\left|\vm{s}_{(i)}\right|_{L_2}$. Then we use (\ref{eq:eAprox}) to evaluate $\left|\vm{e}_{(i)}\right|_{L_2}$ for $i=1,...,N$ whereas the value of $\left|\vm{e}_{(N)}\right|_{L_2}$ is evaluated by extrapolation from the value of $\left|\vm{e}_{(N-1)}\right|_{L_2}$.

This procedure was enhanced in \cite{rodenasgonzalez2007} by evaluating the correction terms $ES_{(i)}$ for $i=1,...,N$ as
\begin{equation}
ES_{(i)}=-2\int_\Omega \vm{e}_{(i)}\cdot\vm{s}_{(i)} d\Omega \approx -2\int_\Omega \left(\vm{u}^h_{(N)}-\vm{u}^h_{(i)} \right)\cdot\vm{s}_{(i)} d\Omega, \ i=1,...,N-1
\label{eq:es-withoutCS}
\end{equation}

and then $ES_{(N)}$ by extrapolation from $ES_{(N-1)}$. 

The two versions of this first alternative have a high computational cost as \textit{i}) it requires the use of a sequence of meshes and \textit{ii})it requires projecting the information from mesh $i$ to mesh $N$. Note also that the results obtained for the first meshes are quite accurate but this accuracy decreases in the last meshes, especially in the last mesh where the results have been obtained using extrapolation techniques. 


The second alternative consists in using a recovery technique that would directly provide an enhanced displacements field $\vm{u}^*$ that could then be used  to estimate the error in the displacements field as $\be \approx \bu^*-\bu^h$ without requiring a sequence of meshes. In this section we will describe the key ideas of the recovery technique used to obtain the error estimates. The technique is also based on the SPR scheme \cite{zienkiewiczzhu1992}, but considering the displacements field and the local satisfaction of the imposed displacement, boundary tractions and internal equilibrium.

We define each of the components of the recovered displacements on the support (patch) of a vertex node $i$ as a polynomial surface expressed as $u^{*}_{i}=\vm{pa}_i$ where $\vm{p}$ is a polynomial basis and $\vm{a}_i$ are the unknown coefficients. The degree of $\vm{p}$ is one order higher than the degree of the shape functions used for the interpolation of displacements. We use the Langrange Multipliers technique to force $\vm{u}^*_{i}$ to satisfy the Dirichlet boundary conditions if the patch contains the Dirichlet boundary. Furthermore, using Lagrange Multipliers over the coefficients $\vm{a}$ it is possible to enforce restrictions on $\vm{\upsigma}^{*}_{i}(x,y)=\vm{DL}\vm{u}^{*}_{i}(x,y)$ so that  $\vm{u}^*_{i}$  satisfies the following equations:
\begin{itemize}
 \item Internal equilibrium equation.
 \item Boundary equilibrium equation.
 \item Compatibility equation (\textit{by default}).
\end{itemize}


The displacements in each element are evaluated from each of the vertex nodes of the element. To enforce continuity of the recovered displacements  we use a partition of unity approach, defined by Belytschko as the Conjoint Polynomial enhancement \cite{blackerbelytschko1994}:

\begin{equation}
\vm{u}^{*}_{u}(x,y)=\sum^{n_v}_{i=1}N_{i}^v(x,y)\vm{u}^{*}_{i}(x,y)
\label{eq:ConjointU}
\end{equation}

Where $N_{i}^v$ represent the shape function of node $i$ corresponding to the linear version of the element (note that patches are only formed around vertex nodes) and $n_v$ refers to the number of vertex nodes. Equation \eqref{eq:ConjointU} provides a kinematically 
admissible displacement field $\vm{u}^{*}_{u}$ evaluated as a function of local contributions $\vm{u}^{*}_{i}$ from each of the different patches. We will use subindex $u$ to indicate a  kinematically admissible solution whereas subindex $\sigma$ will be used to indicate a solution with a continuous stress field. 


With \eqref{eq:StressU}, we can evaluate the stresses $\vm{\upsigma}^{*}_{u}$ at each point as 

\begin{equation}
\begin{split}
\vm{L}&=\begin{bmatrix}
\frac{\partial}{\partial x} & 0 \\
0 & \frac{\partial}{\partial y}  \\
\frac{\partial}{\partial y} & \frac{\partial}{\partial x}
\end{bmatrix} \\
\vm{\upsigma}^{*}_{u}(x,y)&=\vm{DL}\sum^{n_v}_{i=1}\vm{N}_{i}^v(x,y)\vm{u}^{*}_
{ i } (x,y)\\
\vm{\upsigma}^{*}_{u}(x,y)&=\underbrace{\sum^{n_v}_{i=1}\vm{DL}\vm{N}_{i}^v(x,
y)\vm{u}^{*}_{i}(x,y)}_{discontinuous} +
\underbrace{\sum^{n_v}_{i=1}\vm{N}_{i}^v(x,y)\vm{DL}\vm{u}^{*}_{i}(x,y)}_{
\vm{\upsigma}^{*}_{\upsigma}(x,y) \; continuous}  \\
\end{split}
\label{eq:StressU}
\end{equation}

The expression to evaluate $\vm{\upsigma}^{*}_{u}$ has two parts, one discontinuous and another discontinuous. Taking only the continuous part we will have a continuous stress field $\vm{\upsigma}^{*}_{\upsigma}$ which is nearly-statically admissible as it is built using local contributions where the satisfaction of equilibrium has been considered.

We refer to $\vm{\upsigma}^{*}_{\upsigma}$ as nearly-statically admissible field because due to the enforcement of continuity, a lack of internal equilibrium $\vm{s}$ is introduced as shown in \eqref{eq:IntEqU}

\begin{equation}
\begin{split}
\vm{\upsigma}^{*}_{\upsigma}(x,y)&=
\sum^{n_v}_{i=1}\vm{N}_{i}^v(x,y)\vm{DL}\vm{u}^{*}_{i}(x,y) \\
\nabla\vm{\upsigma}^{*}_{\upsigma}(x,y)&=\nabla\sum^{n_v}_{i=1}\vm{N}_{i}^v(x,y)\vm{DL}\vm{u}^{*}_{i}(x,y) \\
\nabla\vm{\upsigma}^{*}_{\upsigma}(x,y)&=\underbrace{\sum^{n_v}_{i=1}
\nabla\vm {N}_{i}^v (x,y)\vm{DL}\vm{u}^{*}_{i}(x,y)}_{-\vm{s}} +
\sum^{n_v}_{i=1}\vm{N}_{i}^v(x,y)\underbrace{\nabla\vm{DL}\vm{u}^{*}_{i}(x,y)}_
{\nabla \vm{\upsigma}^{*}_{i}=-\vm{b}} \\
\nabla\vm{\upsigma}^{*}_{\upsigma}(x,y)&=-\vm{s} - \vm{b} \\
\end{split}
\label{eq:IntEqU}
\end{equation}

Moreover, there could also be a lack of contour equilibrium $\vm{r}$ over $\Gamma_{N}$, $\vm{r} = \vm{\upsigma}^{*}_{\upsigma} \cdot \vm{n} - \vm{t}$, where $\vm{n}$ is the outward normal vector, if the polynomial used to describe the recovered field is not able to represent the applied boundary tractions $\vm{t}$ and/or if the boundary is not a straight line segment.
%
%
\section{Error norm representation for the recovered solution}
\label{sec:ErrRecovered}
%
%
At the same time we  worked on an expression which allow us to evaluate the error of the nearly-statically  recovered stress field $\sst_\upsigma$. Recall that this is a continuous stress field but not fully equilibrated, where $\vm{s}$ and $\vm{r}$ represent the lack of internal and boundary equilibrium, respectively.

The recovery procedure yields a postprocessed stress field $\sst_{\upsigma}$ which is taken as an enhanced approximation to the exact stresses,
$\bsigma(\bu)$, much more accurate than $\bsigma(\bu^h)$. The
recovered stress $\sst_{\upsigma}$ has an equilibrium default 
\begin{equation}
 \vm{s} = - \nabla \cdot \sst_{\upsigma} - \vm{b}
\end{equation}

Pedro Díez proved that the following  identity holds to evaluate the error in the recovered field:
\begin{equation}\label{Eq:StError}
  E^* = -\int _{\Omega} \vm{s} \cdot \vm{e}^{*}_{\upsigma} d\Omega = \norm{\be^{*}_{\upsigma}}^2=
\norm{\bu - \bu^{*}_{\upsigma}}^2 =\bar{a} \left((\bsigma - \sst_{\upsigma}),(\bsigma - \sst_{\upsigma}) \right)
\end{equation}


being $\be^{*}_{\upsigma}$ the error in displacements corresponding to the recovered solution $\sst_{\upsigma}$ such that $\sst_{\upsigma}= \bsigma(\bu^{*}_{\upsigma})
=\bsigma(\bu-\be^{*}_{\upsigma}) $. Note that both
$\bu^{*}_{\upsigma}$ and $\be^{*}_{\upsigma}$ are not explicitly evaluated.

\noindent\textbf{Proof.} 

\begin{align}
\label{eq:ProofSE}
-\int _{\Omega} \vm{s} \cdot \be^{*}_{\upsigma} d\Omega =& 
\int _{\Omega} (\nabla \cdot \sst_{\upsigma} + \vm{b}) \cdot \vm{e}^*_{\upsigma} d\Omega \nonumber\\
=& \int _{\Omega} (\nabla \cdot \bsigma(\bu^{*}_{\upsigma}) - \nabla \cdot
\bsigma(\bu)) \cdot \be^{*}_{\upsigma} d\Omega \nonumber\\
=&  \int _{\Omega} (-\nabla \cdot \bsigma(\be^{*}_{\upsigma})) \cdot \be^{*}_{\upsigma}
d\Omega \nonumber\\
=&  \int _{\Omega} \bsigma(\be^{*}_{\upsigma}): \varepsilon(\be^{*}_{\upsigma}) d\Omega -
 \int _{\partial \Omega} \be^{*}_{\upsigma} \cdot  \bsigma(\be^{*}_{\upsigma}) \cdot \vm{n}
d\Gamma \nonumber\\
=& \int _{\Omega} \bsigma(\be^{*}_{\upsigma}): \varepsilon(\be^{*}_{\upsigma}) d\Omega \nonumber\\
=& \norm{\be^{*}_{\upsigma}}^2
\end{align}

Note that it has been assumed that 
\begin{equation}
 \int _{\partial \Omega} \be^{*}_{\upsigma} \cdot  \vm{\sigma}(\be^{*}_{\upsigma}) \cdot \vm{n}
d\Gamma = 0
\label{eq:boundaryterm-ST}
\end{equation}

If we also consider this term we would have:

\begin{equation}
\label{eq:StErrorFull}
\norm{\be^{*}_{\upsigma}}^2 = -\int _{\Omega} \vm{s} \cdot \be^{*}_{\upsigma} d\Omega  + \int _{\partial \Omega} \be^{*}_{\upsigma} \cdot  \bsigma(\be^{*}_{\upsigma}) \cdot \vm{n} d\Gamma
\end{equation}

Considering  $\bsigma(\be^{*}_{\upsigma}) \cdot \vm{n} =\bsigma(\vm{u}-\vm{u}^{*}_{\upsigma}) \cdot \vm{n} = \bsigma(\vm{u}) \cdot \vm{n} -\vm{\sigma}(\vm{u}^{*}_{\upsigma}) \cdot \vm{n} = \vm{t}-\vm{\sigma}^* \cdot \vm{n}$ and then considering (\ref{eq:Neumann}) we have $\bsigma(\be^{*}_{\upsigma}) = -\vm{r}$. Therefore $\norm{\be^{*}_{\upsigma}}^2$ can be evaluated as 

\begin{equation}
\label{eq:StErrorFull}
E^*=\norm{\be^{*}_{\upsigma}}^2 = -\int _{\Omega} \vm{s} \cdot \be^{*}_{\upsigma} d\Omega  - \int _{\partial \Omega} \vm{r} \cdot \be^{*}_{\upsigma}  d\Gamma
\end{equation}

%


Operating with \eqref{Eq:StError} considering C-S inequality we will have

\begin{equation}
		\left\|\be^{*}_{\sigma}\right\|^2 =-\int_\Omega \be^{*}_{\upsigma}\vm{s}\, d\Omega 
		\leq \left|\be^{*}_{\upsigma}\right|_{L_2}\left|\vm{s}\right|_{L_2}
	\label{eq:RecUB-1}
\end{equation}

As $\vm{u}^*_{\sigma}$ would be a recovered displacements field, laying in a so-called 'broken space' richer than $V^h$, we could assume that the $L_2$ norm of the recovered solution is smaller than the $L_2$ norm of the FE solution:

\begin{equation}
	\left|\be^{*}_{\upsigma}\right|_{L_2} \lesssim \left|\be \right|_{L_2}
	\label{eq:RecUB-2}
\end{equation}

We would then obtain the following practical expression to evaluate an upper bound of $\left\|\be^{*}_{\upsigma}\right\|^2$ taking into account that the error of the FE solution can be estimated as $\vm{e} \approx \vm{e}_{es}$ with $\vm{e}_{es}:=\vm{u}^*_u - \vm{u}^h$

\begin{equation}
	\left\|\be^{*}_{\upsigma}\right\|^2 \leq \left|\be^{*}_{\upsigma}\right|_{L_2}\left|\vm{s}\right|_{L_2} \lesssim \left|\be \right|_{L_2}\left|\vm{s}\right|_{L_2} \approx  \left|\vm{u}^*_u-\vm{u}^h\right|_{L_2}\left|\vm{s}\right|_{L_2} = 
\left|\vm{e}_{es}\right|_{L_2}\left|\vm{s}\right|_{L_2} =E^{*}_{UB}
	\label{eq:RecUB-3}
\end{equation}

%

Note that in (\ref{eq:RecUB-3}) we have finally replaced $\left|\be^{*}_{\upsigma}\right|_{L_2}$ by $\left|\vm{e}_{es}\right|_{L_2}$, \textit{i.e.} we have replaced the error in the recovered solution $\vm{u}^*_{\sigma}$ by the estimated error of the FE solution $\vm{u}^{h}$, to obtain a bound of the error in the recovered solution. Following this idea  we replace $\be^{*}_{\upsigma}=\bu-\bu^{*}_{\upsigma}$ by $\vm{e}_{es}=\vm{u}^*_{u}-\vm{u}_{h}$ in (\ref{eq:StErrorFull}), and defined the error indicator $E^*_1$ in (\ref{eq:RecEst1}) to check if it could provide an indication of the error level in energy norm of the recovered solution $\vm{\sigma}^*_{\upsigma}$. We also defined the error indicators $E^*_2$ and $E^*_3$ as described in (\ref{eq:RecEst2}) and (\ref{eq:RecEst3}) to force the result to be positive:

\begin{subequations} \label{eq:RecEst}
\begin{align}
  E^*_1 &= -\int _{\Omega} \vm{s} \cdot \vm{e}_{es} d\Omega - \int _{\Gamma} \vm{r} \cdot \vm{e}_{es} d\Gamma \label{eq:RecEst1}\\
  E^*_2 &= \sum_{i=1}^{Nel} \left( \abs{\int _{\Omega_i} \vm{s} \cdot \vm{e}_{es} d\Omega} + \abs{\int _{\Gamma\cap\Omega_i} \vm{r} \cdot \vm{e}_{es} d\Gamma} \right) \label{eq:RecEst2}\\
  E^*_3 &= \int _{\Omega} \abs{\vm{s} \cdot \vm{e}_{es}}  d\Omega + \int _{\Gamma} \abs{\vm{r} \cdot \vm{e}_{es}} d\Gamma \label{eq:RecEst3}
\end{align}
\end{subequations}


where $Nel$ is the number of elements of the FE discretization and $\Omega_i$ represents the domain of element $i$.

In (\ref{eq:RecEst2}) the value of the integrals at each element are forced to be positive. In (\ref{eq:RecEst3}) the integrands themselves are forced to be positive. Note that this is a reasonable assumption.  For example, if we assume $\vm{r}=\vm{0}$ (see (\ref{eq:ProofSE}))

\begin {equation}
0\leq\bsigma(\be^{*}_{\upsigma})\varepsilon(\be^{*}_{\upsigma}) = -\vm{s}\cdot\be^{*}_{\upsigma}\ \  \left(\:= \left| \vm{s} \cdot \be^{*}_{\upsigma}\right| \right)
\end {equation}

As $\vm{s}$ and $\be^{*}_{\upsigma}$ are consistent ($\vm{s}$ would be the defaults of equilibrium corresponding to $\vm{u}_{\upsigma}^*$, whose associated error is $\be^{*}_{\upsigma}$) then $0 \le -\vm{s} \cdot \be^{*}_{\upsigma}$. However, in (\ref{eq:RecEst3}) $\be^{*}_{\upsigma}$ has been substituted by $\vm{e}_{es}$. The terms $\vm{s}$ and $\vm{e}_{es}$ are non-consistent and as a result $-\vm{s} \cdot \vm{e}_{es}$ could be negative. This suggested the use of the approximation in (\ref{eq:RecEst3}) ($-\vm{s} \cdot \be^{*}_{\upsigma}\approx\left| \vm{s} \cdot \be_{es}\right|$).

We wanted to check if the error indicators $E^*_i$'s could provide a rough idea about the error level of the recovered solution. We also wanted to check if the local evaluations, at element level, of the $E^*_i$ could roughly describe the distribution of the error of the recovered solution. As far as we understand, there is no mathematical relation between $E^*$ and the error indicators $E^*_i$ that could make them provide accurate evaluations of $E^*$. However, the numerical results obtained show that the $E^*_i$'s capture the order of magnitude of $E^*$. In particular, $E^*_3$ provides an accurate evaluation of $E^*$ ( $E^*_3 \approx E^*$) and very similar error distributions.
%
%
\section{Numerical Examples}
\label{sec:Examples}
%
%

Quadrilateral elements have been considered in the examples presented in this section.
%
%
\subsection{Westergaard problem. XFEM solution}
\begin{figure}[h!]
    \centering
    \includegraphics{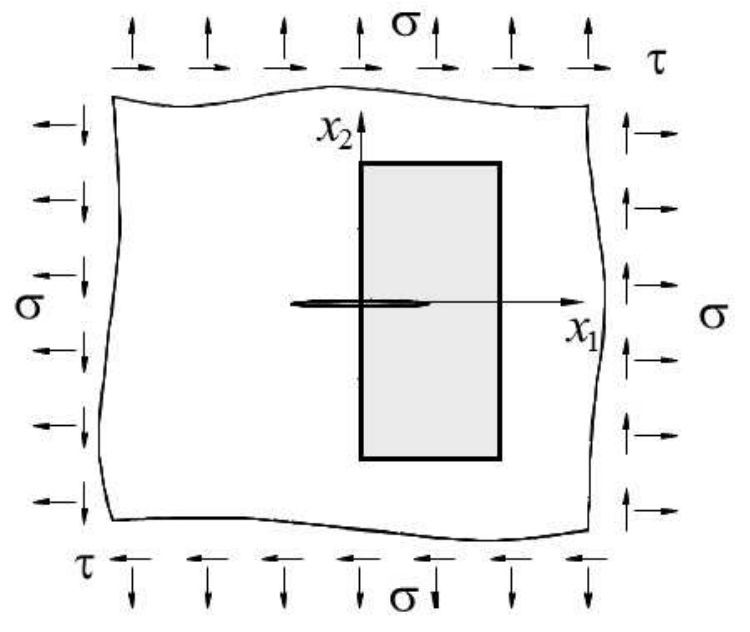}
    \caption{Example 1. Westergaard problem. Infinite plate with a crack of length $2a$ under uniform tractions $\sigma_{\infty}$ (biaxial) and $\tau_{\infty}$. Finite portion of the domain $\Omega_0$, modelled with FE.}
    \label{fig:westergaard}
\end{figure}

%
%
Our first results were obtained as part of the development of an error estimator for XFEM. Let us consider the Westergaard problem corresponding to an infinite plate loaded at infinity with biaxial tractions $\sigma_{x \infty}=\sigma_{y \infty}=\sigma_{\infty}$ and shear traction $\tau_{\infty}$, presenting a crack of length $2a$ as shown in Figure~\ref{fig:westergaard}. The problem is solved using an enriched finite element approximation, combining the externally applied loads to obtain different loading conditions: pure mode I, pure mode II or mixed mode.  
 
The numerical model corresponds to a finite portion of the domain ($a=1$ and $b=4$ in Figure~\ref{fig:westergaard}). The applied projected stresses for mode I are evaluated from the analytical Westergaard solution \cite{ginerfuenmayor2005}:  

\begin{equation} \label{Eq:WesterStressI} 
\begin{array}{r@{\hspace{1ex}}c@{\hspace{1ex}}l} 
{\sigma _{x}^{I} } (x,y) & {=} & {\displaystyle \frac{\sigma _{\infty } }{\sqrt{\left|t\right|} } \bigg[\left(x\cos \frac{\phi }{2} -y\sin \frac{\phi }{2} \right)+y\frac{a^{2} }{\left|t\right|^{2} } \left(m\sin \frac{\phi }{2} -n\cos \frac{\phi }{2} \right)\bigg]} \\ 
\noalign{\medskip}{\sigma _{y}^{I} }(x,y) & {=} & {\displaystyle \frac{\sigma _{\infty } }{\sqrt{\left|t\right|} } \bigg[\left(x\cos \frac{\phi }{2} -y\sin \frac{\phi }{2} \right)-y\frac{a^{2} }{\left|t\right|^{2} } \left(m\sin \frac{\phi }{2} -n\cos \frac{\phi }{2} \right)\bigg]} \\ 
\noalign{\medskip}{ \tau _{xy}^{I} }(x,y) & {=} & {\displaystyle y\frac{a^{2} \sigma _{\infty } }{\left|t\right|^{2} \sqrt{\left|t\right|} } \left(m\cos \frac{\phi }{2} +n\sin \frac{\phi }{2} \right)} \end{array}
\end{equation}
\noindent and for mode II:
\begin{equation} \label{Eq:WesterStressII} 
\begin{array}{r@{\hspace{1ex}}c@{\hspace{1ex}}l} 
{\sigma _{x}^{II}}(x,y)  & {=} & {\displaystyle \frac{\tau _{\infty } }{\sqrt{\left|t\right|} } \bigg[2\left(y\cos \frac{\phi }{2} +x\sin \frac{\phi }{2} \right)-y\frac{a^{2} }{\left|t\right|^{2} } \left(m\cos \frac{\phi }{2} +n\sin \frac{\phi }{2} \right)\bigg]} \\ 
\noalign{\medskip}{\sigma _{y}^{II}}(x,y)  & {=} & {\displaystyle y\frac{a^{2} \tau _{\infty } }{\left|t\right|^{2} \sqrt{\left|t\right|} } \left(m\cos \frac{\phi }{2} +n\sin \frac{\phi }{2} \right)} \\ \noalign{\medskip}{\tau _{xy}^{II}}(x,y)  & {=} & {\displaystyle \frac{\tau _{\infty } }{\sqrt{\left|t\right|} } \bigg[\left(x\cos \frac{\phi }{2} -y\sin \frac{\phi }{2} \right)+y\frac{a^{2} }{\left|t\right|^{2} } \left(m\sin \frac{\phi }{2} -n\cos \frac{\phi }{2} \right)\bigg]} \end{array}
\end{equation}

\noindent where the stress fields are expressed as a function of $x$ and $y$, with origin at the centre of the crack. The parameters $t$, $m$, $n$ and $\phi$ are defined as
\begin{equation}
\begin{split}t& =(x+iy)^{2} -a^{2} =(x^{2} -y^{2} -a^{2} )+i(2xy)=m+in \\  m & =\textrm{Re}(t) =\textrm{Re}(z^{2} -a^{2} )=x^{2} -y^{2} -a^{2} \\ n & =\textrm{Im}(t)=(z^{2} -a^{2} )=2xy \\ \phi & =\textrm{Arg} (\bar{t})=\textrm{Arg} (m-in) \qquad\textrm{with }\phi \in \left[-\pi ,\pi \right], \; i^2=-1 \end{split} 
\end{equation}

For the problem analysed, the exact value of the SIF is given by 
\begin{equation} \label{Eq:SIFWestergaard} 
K_{I,ex} =\sigma_{\infty } \sqrt{\pi a} \qquad \qquad K_{II,ex} =\tau_{\infty } \sqrt{\pi a}  
\end{equation}

Material parameters are Young's modulus $E = 10^7$ and Poisson's ratio $ \nu= 0.333$. We consider loading conditions in pure mode I with $\sigma_{\infty} =100$ and $\tau_{\infty}=0$, pure  mode II with $\sigma_{\infty} =0$ and $\tau_{\infty}=100$, and mixed mode with $\sigma_{\infty} =100$ and $\tau_{\infty}=100$. In the numerical analyses, we use  a geometrical enrichment defined by a circular fixed enrichment area $B(x_0, r_e)$ with radius  $r_e = 0.5$, with its centre at the crack tip $x_0$ as proposed in \cite{bechetminnebo2005}. For the extraction of the SIF we define a plateau function with radius $r_q = 0.9$. Bilinear elements are considered in the models. For the numerical integration of standard elements we use a $2\times2$ Gaussian quadrature rule. The elements intersected by the crack are split into triangular integration subdomains that do not contain the crack. We use 7 Gauss points in each triangular subdomain, and a $5\times5$ quasipolar integration in the subdomains of the element containing the crack tip. We do not consider correction for blending elements.

Note that for this problem solved using XFEM, the displacements recovery technique was not implemented, therefore we used a stress recovery technique called \mbox{SPR-C}, see \cite{diezrodenas2007} (this technique enforces the recovered stress field to satisfy at each patch the internal and boundary equilibrium equations and the compatibility equation. The recovery procedure also splits the stress field into singular and smooth parts:  $\vm{\upsigma}=\vm{\upsigma}_{smo}+\vm{\upsigma}_{sing}$, and then uses different recovery procedures for each part as described in \cite{rodenasgonzalez2008} and \cite{rodenasgonzalez2007}. Therefore, we  used the estimation of the error in the displacements field given in (\ref{eq:eAprox}), and used the following expression (equivalent to equation (\ref{eq:RecUB-3})) to evaluate an upper bound of the error of the recovered solution for a sequence of $N$ meshes

\begin{align}
	\left\|\be^{*}_{(i)\upsigma}\right\|^2 &\leq \left|\be^{*}_{(i)\upsigma}\right|_{L_2}\left|\vm{s}_{(i)}\right|_{L_2} \lesssim \left|\be_{(i)} \right|_{L_2}\left|\vm{s}_{(i)}\right|_{L_2}\approx\left|\vm{u}^h_{(N)}-\vm{u}^h_{(i)}\right|_{L_2}\left|\vm{s}_{(i)}\right|_{L_2}  \nonumber \\
&= \left|\vm{e}_{(i)es}\right|_{L_2}\left|\vm{s}_{(i)}\right|_{L_2} =E^{*}_{(i)UB}\ \ \ \ \ \ \ \ i=1,...,N-1
	\label{eq:RecUB-3XFEM}
\end{align}
 
then, we evaluated $\left|\vm{e}_{(N)es}\right|_{L_2}$ by extrapolation from $\left|\vm{e}_{(i-1)es}\right|_{L_2}$.

The error in energy norm of the recovered solution was evaluated using the following expression (similar to (\ref{eq:RecEst3})) where the term associated to the lack of boundary equilibrium was neglected because $\vm{r}\approx \vm{0}$:

\begin{equation} \label{eq:RecEst3XFEM}
 E^*_{(i)3} = \int_{\Omega}\abs{ \vm{s}_{(i)}\cdot\left( \vm{u}^h_{(N)}-\vm{u}^h_{(i)} \right) } d\Omega = \int _{\Omega} \abs{\vm{s}_{(i)} \cdot \vm{e}_{(i)es}}  d\Omega; \ \ \ i=1,...,N-1, 
\end{equation}

we evaluated $E^*_{(N)3}$ by extrapolation from $E^*_{(N-1)3}$.

\begin{figure}[h!]
\centering
  \includegraphics[width=0.45\textwidth]{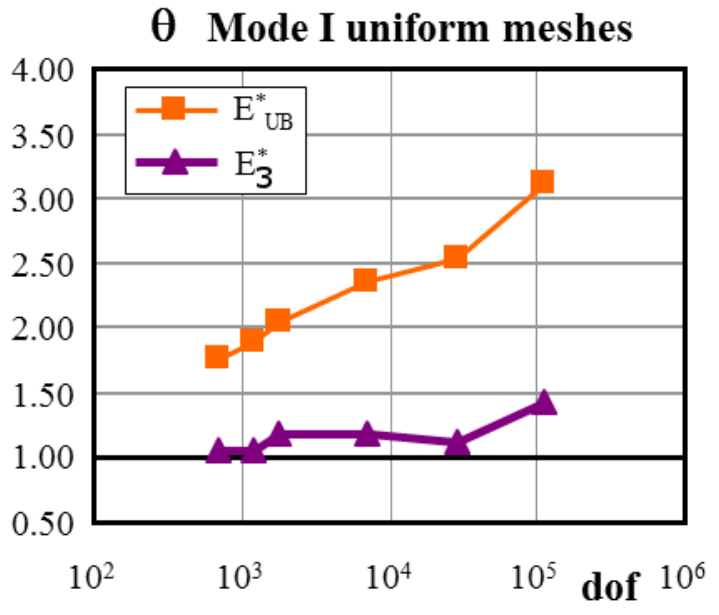} \includegraphics[width=0.45\textwidth]{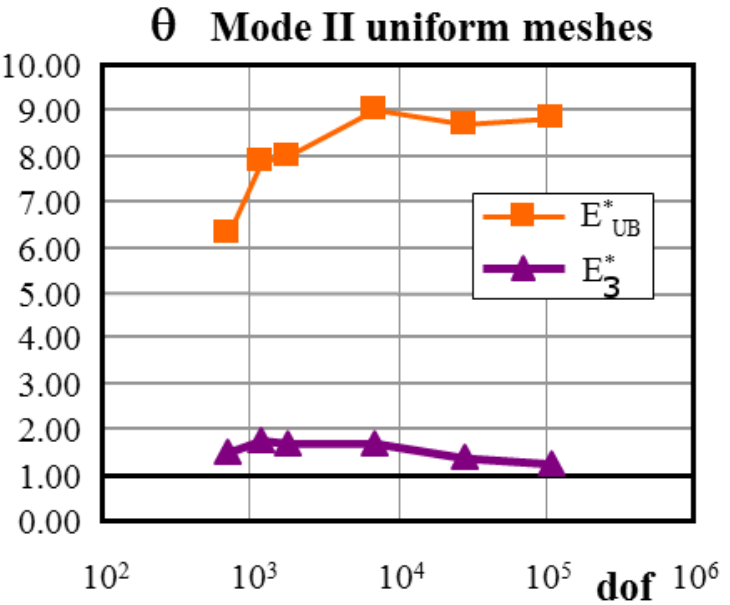}  
\includegraphics[width=0.45\textwidth]{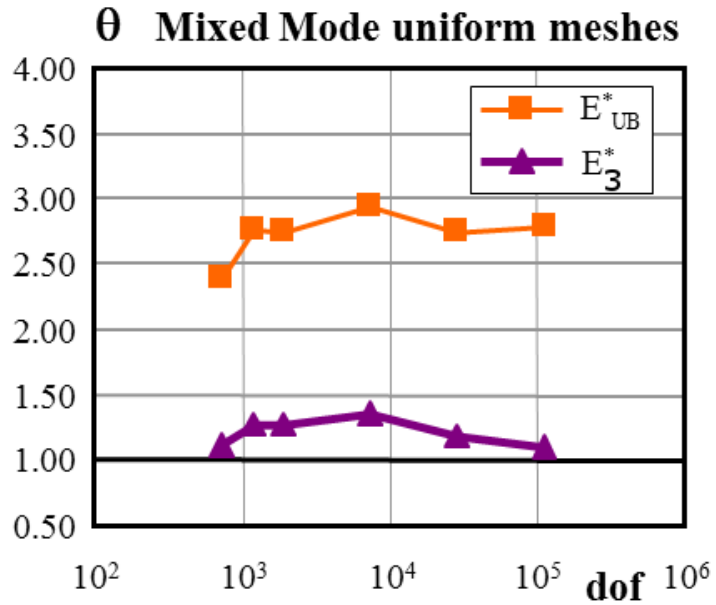}
  \caption{Example 1. Westergaard problem. Effectivity for the estimate of the error of the recovered solution in mode I, mode II and mix mode.}
  \label{fig:ThetaWestergaard}
\end{figure}

Figure \ref{fig:ThetaWestergaard} shows the evolution of the global effectivity index for the error estimate $E_3^*$ and the upper bound $E^*_{\rm UB}$. We can observe that although $E^*_{\rm UB}$ bounds the exact error, it does not asymptotically converge to 1. We believe that we obtain bounds of the exact error because of the use of the CS inequality. Nevertheless, we can observe that the effectivity index of the error estimator in energy norm for the recovered solution $E^*_3$ captures the order of magnitude of the exact error, providing effectivities very close to unity. In some cases, as in the top left graph, we can obtain worse accuracies in the last mesh because of the extrapolation procedure used.

%
%
%
%
%
%
%
%
%
%
%
%
%
\subsection{Square plate with 4$^{th}$ order displacements }
\begin{figure}[h!]
	\begin{tabular}{m{0.45\textwidth} l}
		\includegraphics[width=0.4\textwidth]{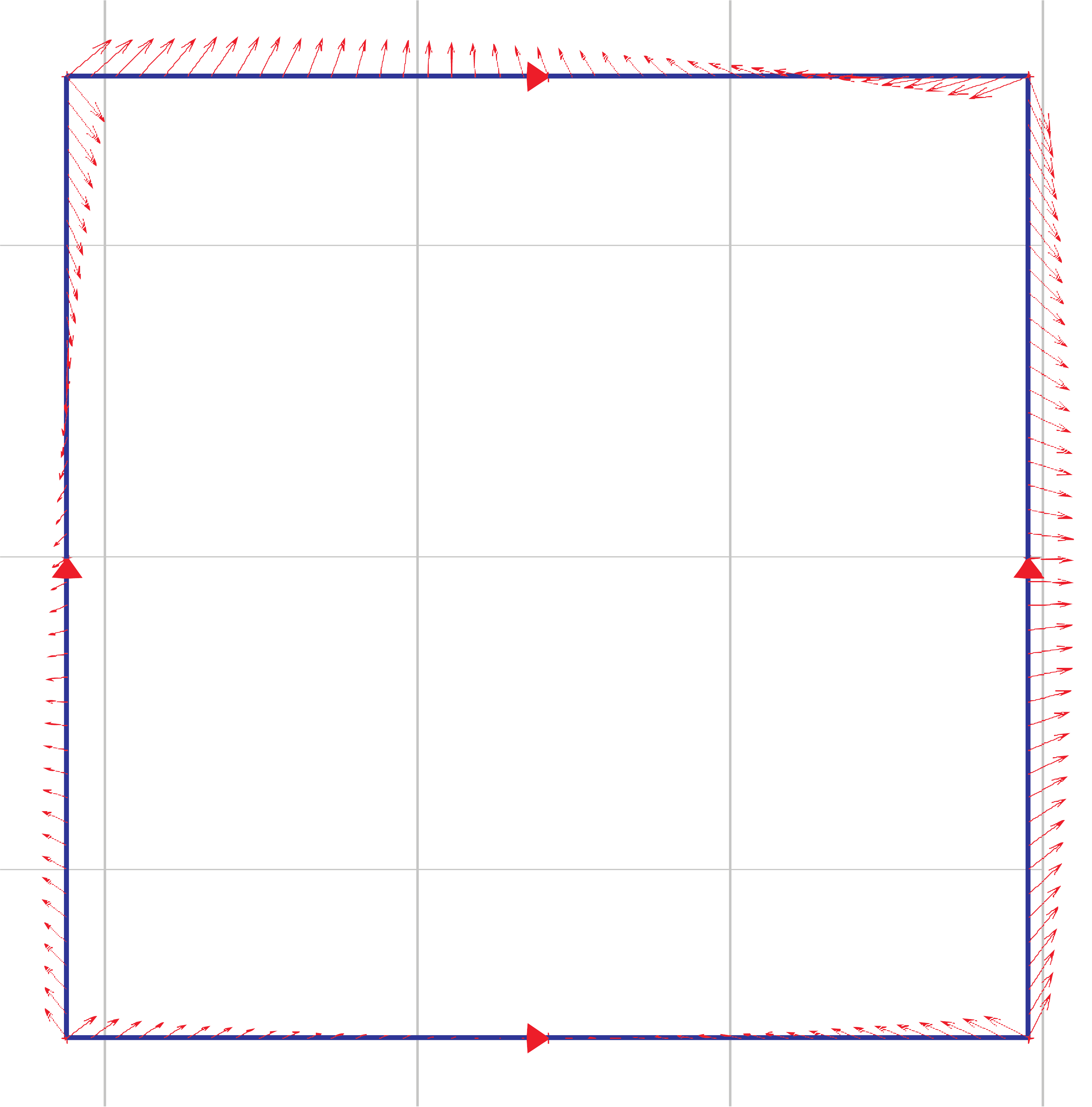} 
	&
		\parbox[l]{.55\textwidth}{
		\begin{align*}
			&u_x(x,y) = x^4+5(x^3y)-3(x^2y^2)+x^3\\
			&u_y(x,y) = y^4-6(y^2x^2)+3(yx^3)+2y\\
			&b_x(x,y) = 	  \frac{-3E}{2(1+\nu)(2\nu-1)}...\\
			&(9x^2-12xy+4y^2-4x+...\\
			&\nu(4x^2+20xy-4y^2+4x))\\
			&b_y(x,y) = \frac{3E}{2(1+\nu)(2\nu-1)}...\\
			&(4y^2-3x^2+2xy+\nu(8x^2-12xy) \\
			&E=1000 \; \nu = 0.3
		\end{align*}
			}
	\end{tabular}
	\caption{Example 2: $2\times2$ square model and analytical solution
		}
	\label{fig:square_model}
\end{figure}

%
%
The FE software used with this and the following problems is based on the use of Cartesian grids independent of the shape of the problem to be analyzed. Elements fully located in the interior of the domain are considered as standard elements. Elements fully outside of the domain are not assembled. Finally, in elements cut by the boundary, element matrices are evaluated considering only the portion of the elements within the domain. The Lagrange multipliers technique has been used to impose a weak satisfaction of the displacement constraints. The  convergence rate of the  error obtained is the theoretical rate expected for the standard FEM. The FE code also uses a hierarchical data structure to perform \textit{h}-adaptive analyses based on element splitting. Multi-point constraints are used to enforce $C^0$ continuity of the solution between adjacent elements of different refinement levels.
This problem considers a 4$^{th}$ order displacements field over an infinite domain. A  2$\times$2 square portion  has been modelled. The corresponding body loads and Neumann conditions have been imposed. A model of the problem is represented in Figure \ref{fig:square_model}.

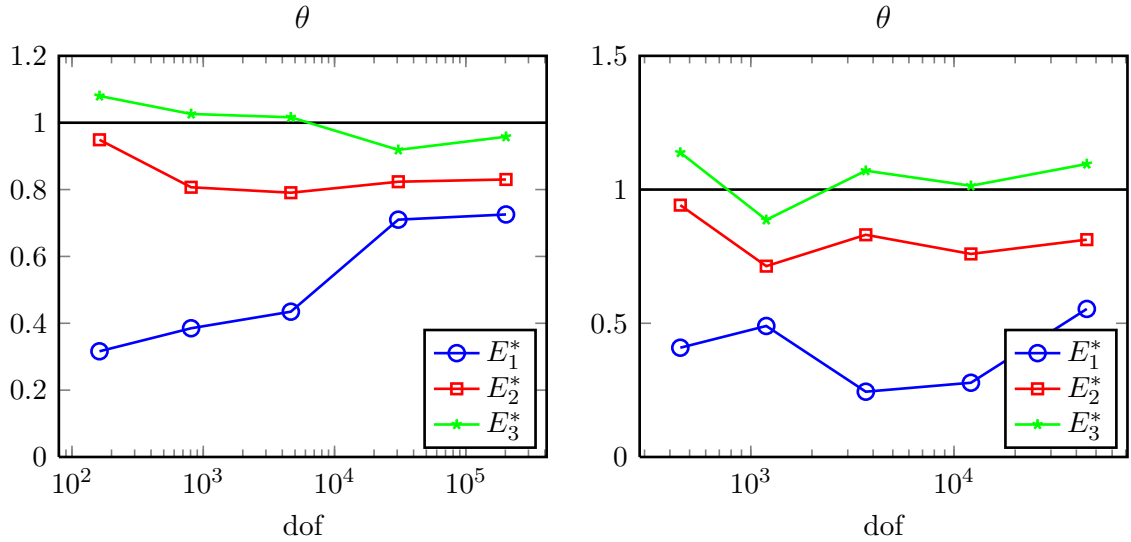
\begin{figure}[h!]
\centering
\begin{tabular}{c c}
   \begin{tikzpicture}
 	\begin{semilogxaxis}[
 	title={$\theta$},
 	ymax=1.2,
 	ymin=0,
 	xlabel={dof},
 	legend style={cells={anchor=west}, font=\small},
 	legend style={at={(0.98,0.02)},anchor=south east},
 	cycle list name=ageplot
 	]
 	\draw (axis cs: 0,1)--(axis cs: 1e8,1); 
 	\addplot table[x=NDOF,y= Eff*_M1]{Data/SQ-4rd-Q4_T005.csv};
 	\addplot table[x=NDOF,y= Eff*_M2]{Data/SQ-4rd-Q4_T005.csv};
 	\addplot table[x=NDOF,y= Eff*_M3]{Data/SQ-4rd-Q4_T005.csv};
 	\legend{{$E^*_1$},{$E^*_2$},{$E^*_3$}}
 	\end{semilogxaxis}
 \end{tikzpicture}
&
   \begin{tikzpicture}
 	\begin{semilogxaxis}[
 	title={$\theta$},
 	ymax=1.5,
 	ymin=0,
 	xlabel={dof},
 	legend style={cells={anchor=west}, font=\small},
 	legend style={at={(0.98,0.02)},anchor=south east},
 	cycle list name=ageplot
 	]
 	\draw (axis cs: 0,1)--(axis cs: 1e8,1); 
 	\addplot table[x=NDOF,y= Eff*_M1]{Data/SQ-4rd-Q8_T006.csv};
 	\addplot table[x=NDOF,y= Eff*_M2]{Data/SQ-4rd-Q8_T006.csv};
 	\addplot table[x=NDOF,y= Eff*_M3]{Data/SQ-4rd-Q8_T006.csv};
 	\legend{{$E^*_1$},{$E^*_2$},{$E^*_3$}}
 	\end{semilogxaxis}
 \end{tikzpicture} 
\end{tabular}
 \caption{Example 1. $2\times2$ square. Effectivity index. Q4 (left) and Q8 (right).} 
 \label{fig:SQ}
\end{figure}

In Figure \ref{fig:SQ} we represent the evolution of the global effectivity index $\theta$  considering the error estimates $E_i^*$ defined in (\ref{eq:RecEst}). Note that in this case we have considered the lack of equilibrium both in the internal and boundary equilibrium equations.

The results show that the error estimators $E^*_1$, $E^*_2$ and $E^*_3$ capture the order of magnitude of the exact error in energy norm of the recovered solution. In particular, the results obtained with $E^*_3$ provide error effectivity indexes very close to one, the desired value. Taking into account the procedure described to obtain these error estimators, the results obtained, especially with $E^*_3$, are surprisingly accurate.

\begin{figure}[!]
  \centering
  \begin{tabular}{c c c}
   Mesh & $\left\|\be^*\right\|^2_k$ & $E^{*}_{3}|_k$ \\
   2 & \includegraphics[trim = 1.0cm 0.4cm 2.1cm 0.3cm, clip, width=0.4\textwidth]{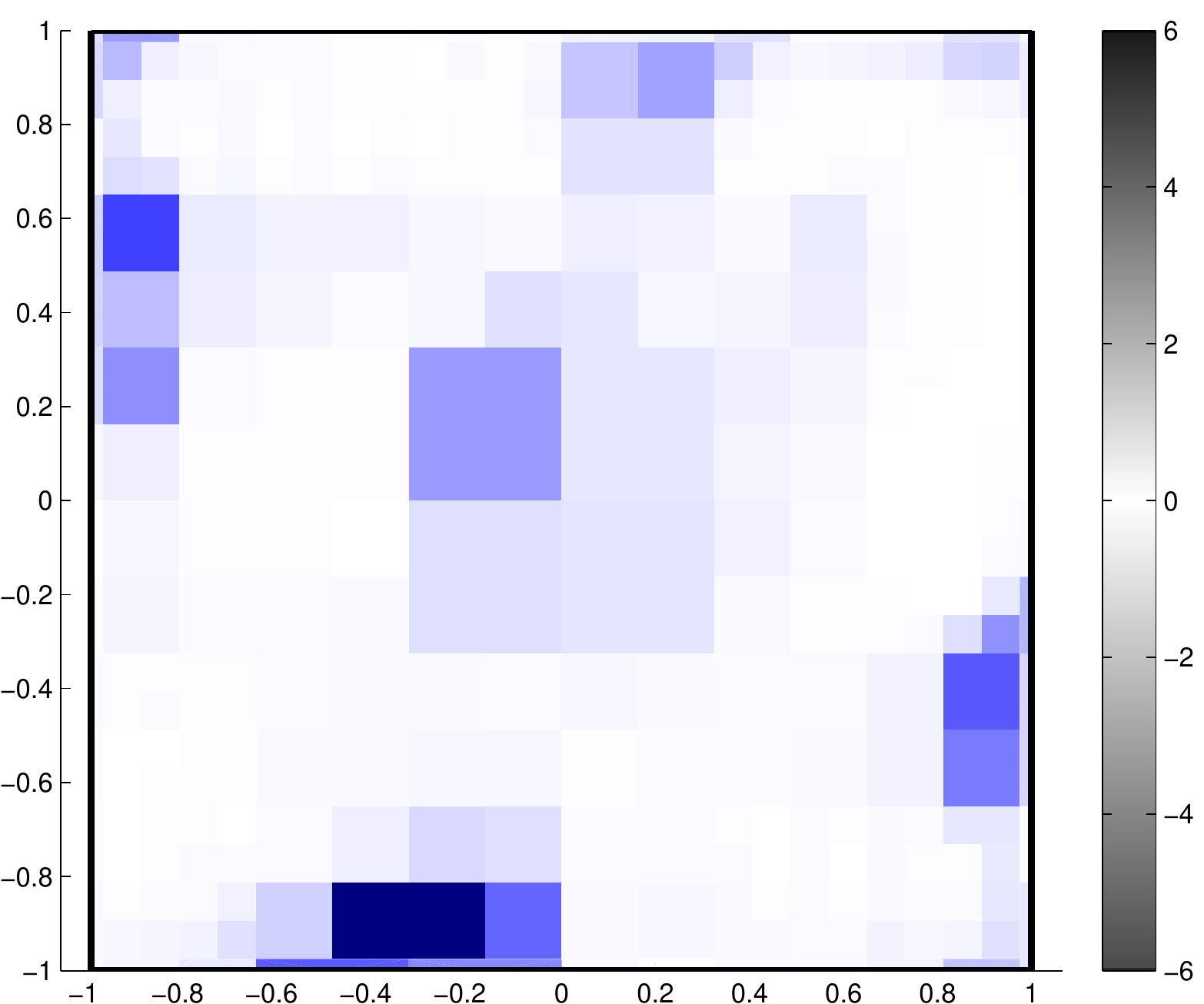} &   	\includegraphics[trim = 1.0cm 0.4cm 2.1cm 0.3cm, clip, width=0.4\textwidth]{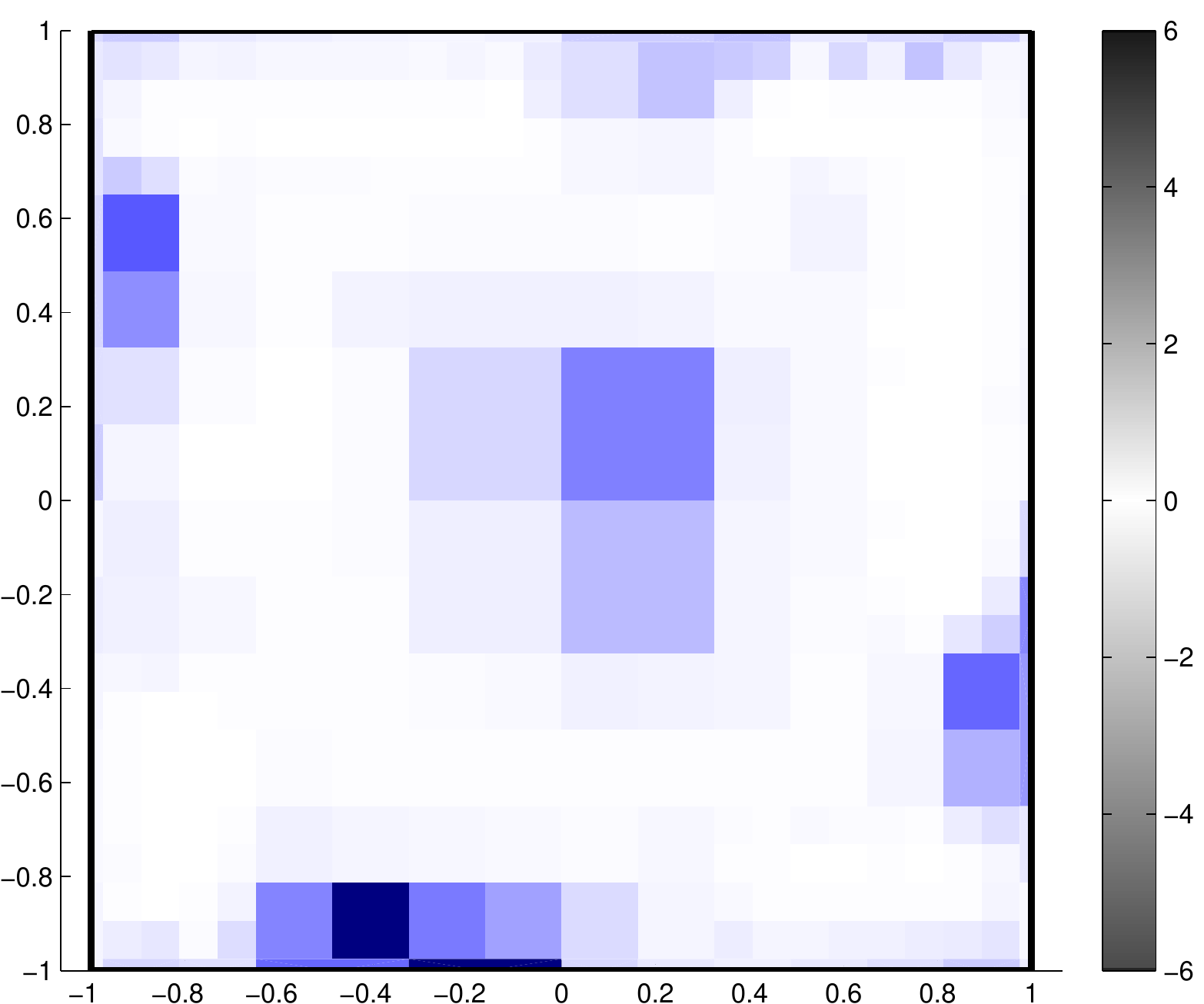} \\
   3 & \includegraphics[trim = 3.2cm 0.4cm 2.6cm 0.3cm, clip, width=0.4\textwidth]{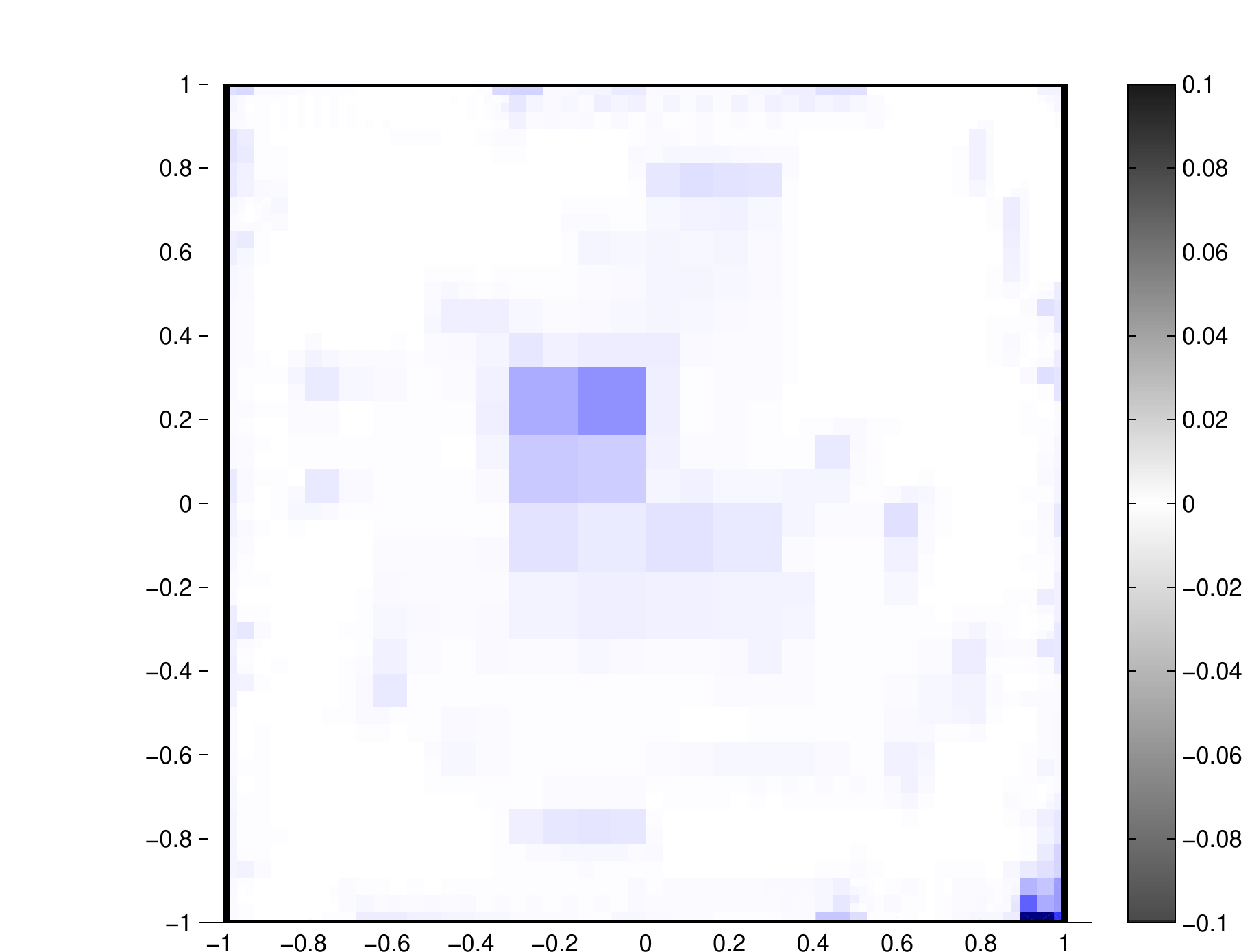} &   	\includegraphics[trim = 3.2cm 0.4cm 2.6cm 0.3cm, clip, width=0.4\textwidth]{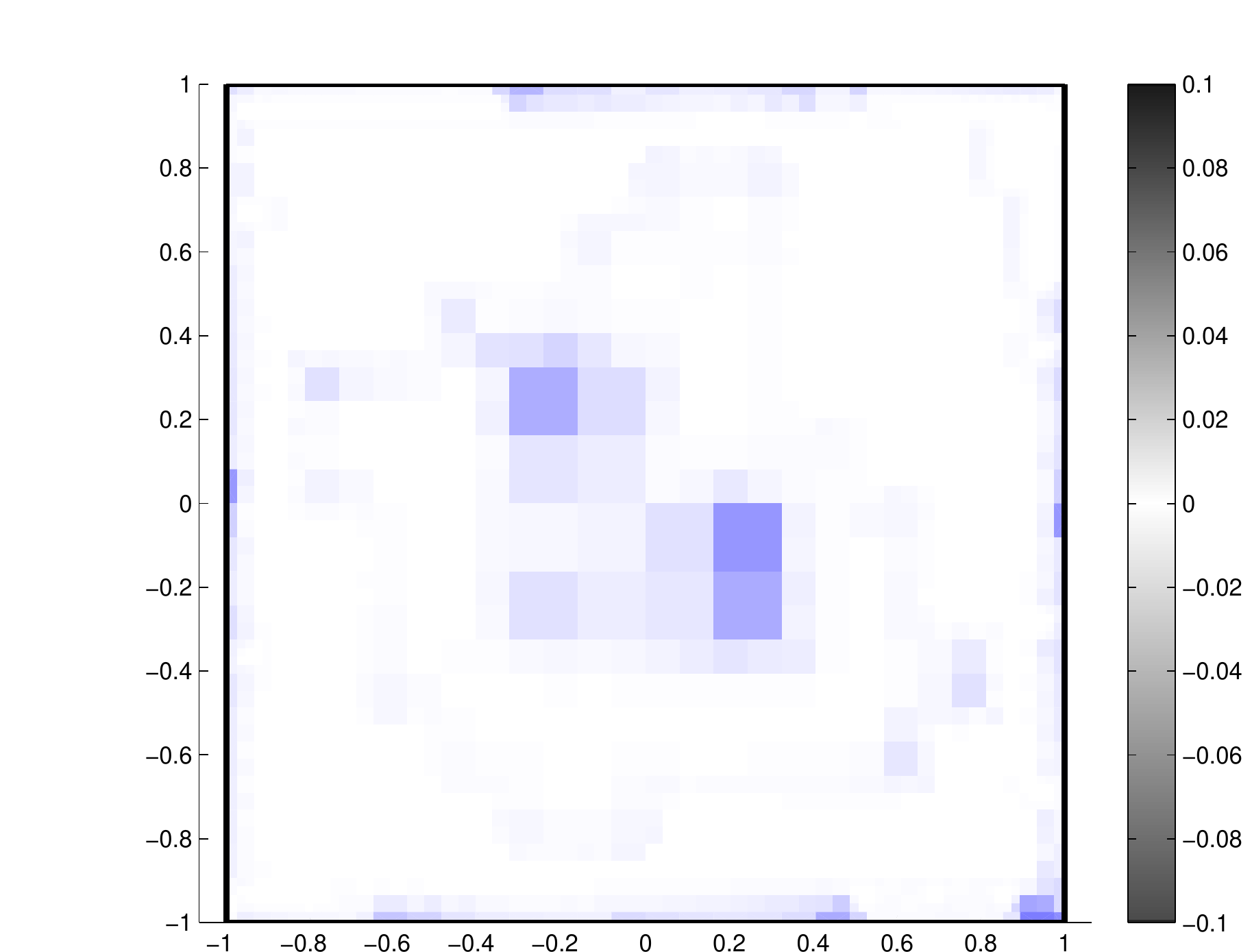} \\
   4 & \includegraphics[trim = 3.2cm 0.4cm 2.6cm 0.3cm, clip, width=0.4\textwidth]{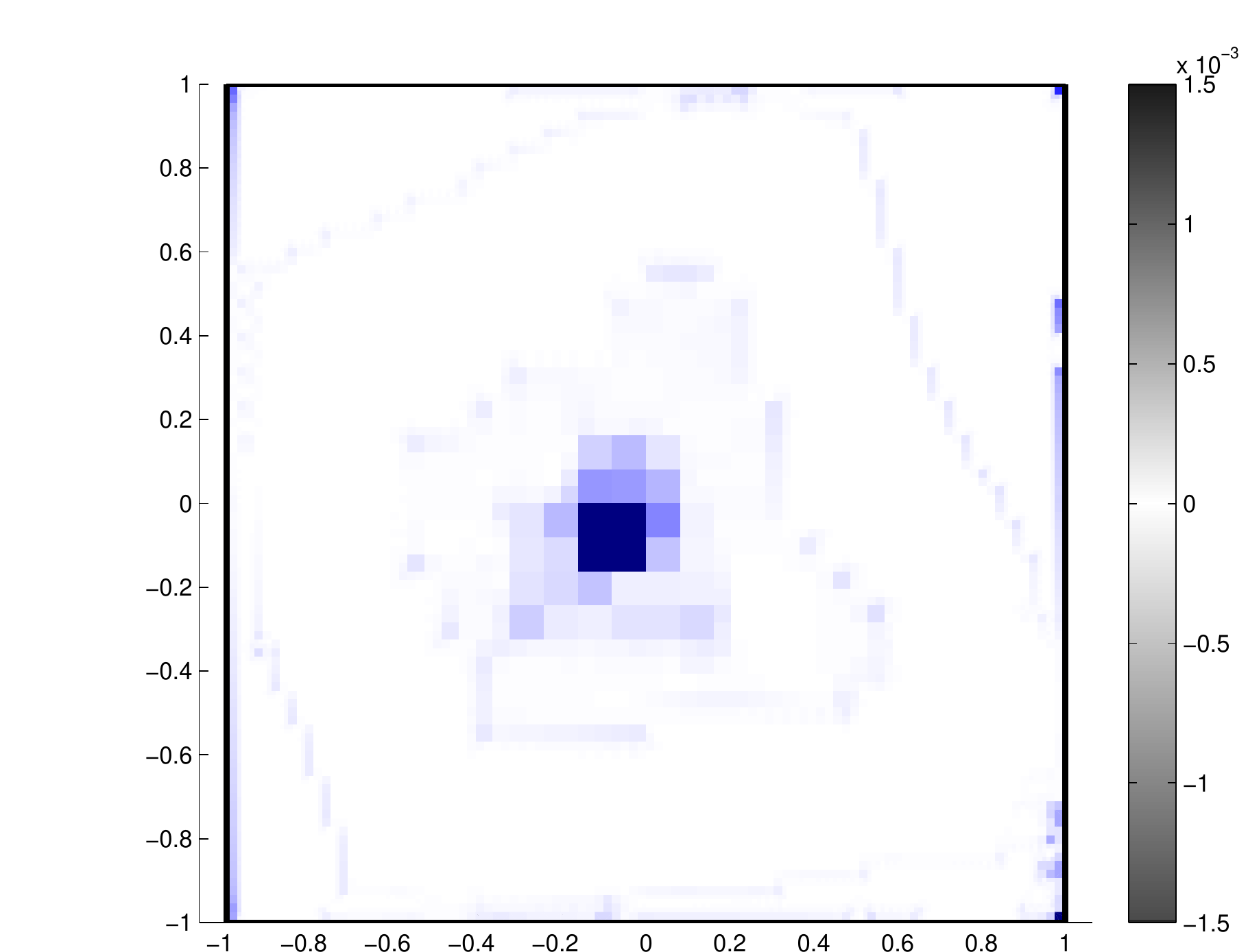} &   	\includegraphics[trim = 3.2cm 0.4cm 2.6cm 0.3cm, clip, width=0.4\textwidth]{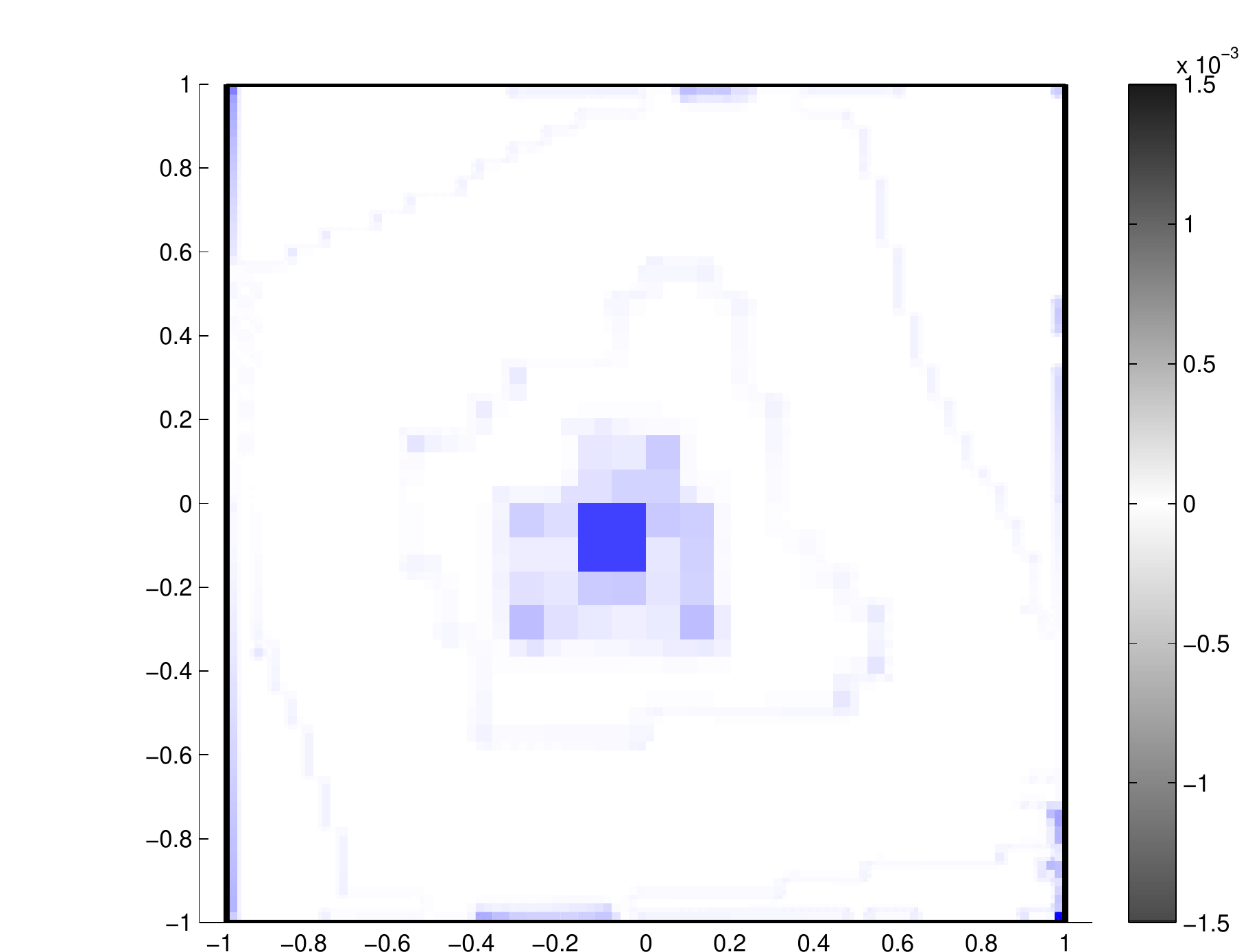} \\
   \end{tabular}
  \caption{Example 2. Square. Q4. Local exact error of the recovered solution (left),  local error estimates using $E^{*}_{3}$ estimator (right).}
  \label{fig:SQLocErrEst}
\end{figure}

Moreover, in Figure \ref{fig:SQLocErrEst} we show the local error evaluated by $E^{*}_{3}|_k$ and the exact error of the recovered field $\left\|\be^*\right\|^2_k$ for a sequence of \mbox{\textit{h}-adapted} meshes. $k$ indicates that the quantities are evaluated at element level. For this problem both results are quite similar, thus $E^{*}_{3}|_k$ is a good indicator at local level of the error of the recovered solution.

%
%
\subsection{Pipe under internal pressure}
\begin{figure}[!]
	\begin{tabular}{m{0.45\textwidth} l}
		\includegraphics[width=0.4\textwidth]{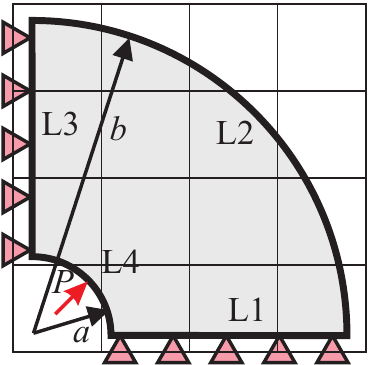} 
		&
		\parbox[l]{.55\textwidth}{
		\begin{align*}
			u_r(r) &=\frac{P(1+\nu)}{E(c^2-1)}\left(r(1-2\nu)+\frac{b^2}{r}\right)\\
			\upsigma_r(r)&=\frac{P}{c^2-1}\left(1-\frac{b^2}{r^2}\right)\\
			\upsigma_\uptheta(r)&=\frac{P}{c^2-1}\left(1+\frac{b^2}{r^2}\right)\\
			a&=5\qquad b=20\qquad P=1 \\ E&=1000\qquad  \nu=0.3 \\
			c&=\frac{b}{a}
		\end{align*}
			}
	\end{tabular}
	\caption{Example 3: Pipe model and analytical solution
		}
	\label{fig:Pipe_model}
\end{figure}

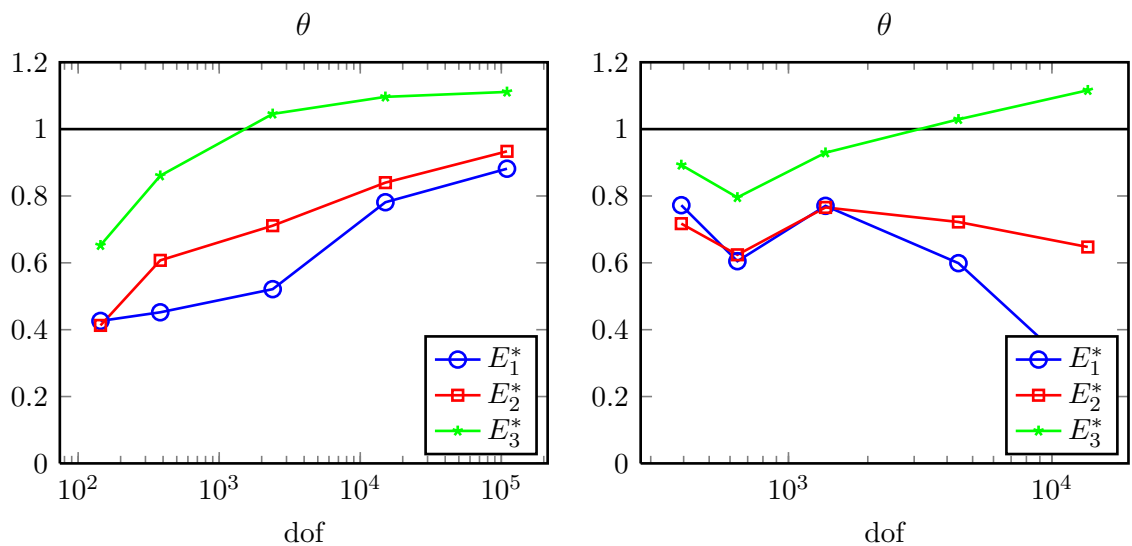
\begin{figure}[!]
\centering
\begin{tabular}{c c}
   \begin{tikzpicture}
 	\begin{semilogxaxis}[
 	title={$\theta$},
 	ymax=1.2,
 	ymin=0,
 	xlabel={dof},
 	legend style={cells={anchor=west}, font=\small},
 	legend style={at={(0.98,0.02)},anchor=south east},
 	cycle list name=ageplot
 	]
 	\draw (axis cs: 0,1)--(axis cs: 1e8,1); 
 	\addplot table[x=NDOF,y= Eff*_M1]{Data/Cyl-Q4_T039.csv};
 	\addplot table[x=NDOF,y= Eff*_M2]{Data/Cyl-Q4_T039.csv};
 	\addplot table[x=NDOF,y= Eff*_M3]{Data/Cyl-Q4_T039.csv};
 	\legend{{$E^*_1$},{$E^*_2$},{$E^*_3$}}
 	\end{semilogxaxis}
 \end{tikzpicture}
&
   \begin{tikzpicture}
 	\begin{semilogxaxis}[
 	title={$\theta$},
 	ymax=1.2,
 	ymin=0,
 	xlabel={dof},
 	legend style={cells={anchor=west}, font=\small},
 	legend style={at={(0.98,0.02)},anchor=south east},
 	cycle list name=ageplot 	]
 	\draw (axis cs: 0,1)--(axis cs: 1e8,1); 
 	\addplot table[x=NDOF,y= Eff*_M1]{Data/Cyl-Q8-T040.csv};
 	\addplot table[x=NDOF,y= Eff*_M2]{Data/Cyl-Q8-T040.csv};
 	\addplot table[x=NDOF,y= Eff*_M3]{Data/Cyl-Q8-T040.csv};
 	\legend{{$E^*_1$},{$E^*_2$},{$E^*_3$}}
 	\end{semilogxaxis}
 \end{tikzpicture} 
\end{tabular}
\caption{Example 3. Pipe. Effectivity index. Q4 (left) and Q8 (right).} 
\label{fig:PipeGolbEfect}
\end{figure}

%
%
Let us now consider a pipe under internal pressure, as described in Figure \ref{fig:Pipe_model}. Figure \ref{fig:PipeGolbEfect} shows the evolution of the effectivity index for the different error estimates. 

\begin{figure}[!]
  \centering
  \begin{tabular}{c c c}
   Mesh & $\left\|\be^*\right\|^2_k$ & $E^{*}_{3}|_k$ \\
   2 & \includegraphics[trim = 0.7cm 0.4cm 2.6cm 0.3cm, clip, width=0.4\textwidth]{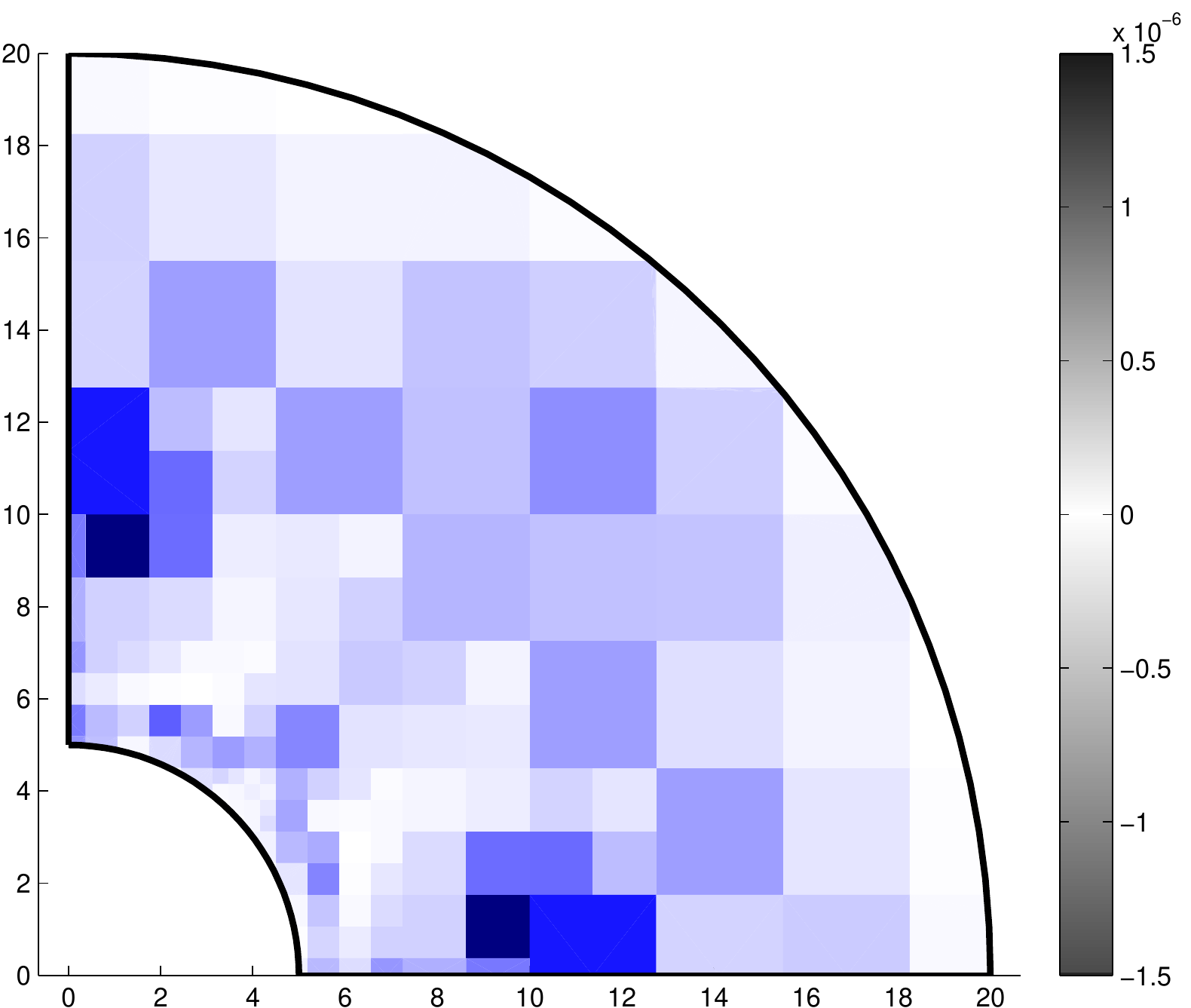} & \includegraphics[trim = 0.7cm 0.4cm 2.6cm 0.3cm, clip, width=0.4\textwidth]{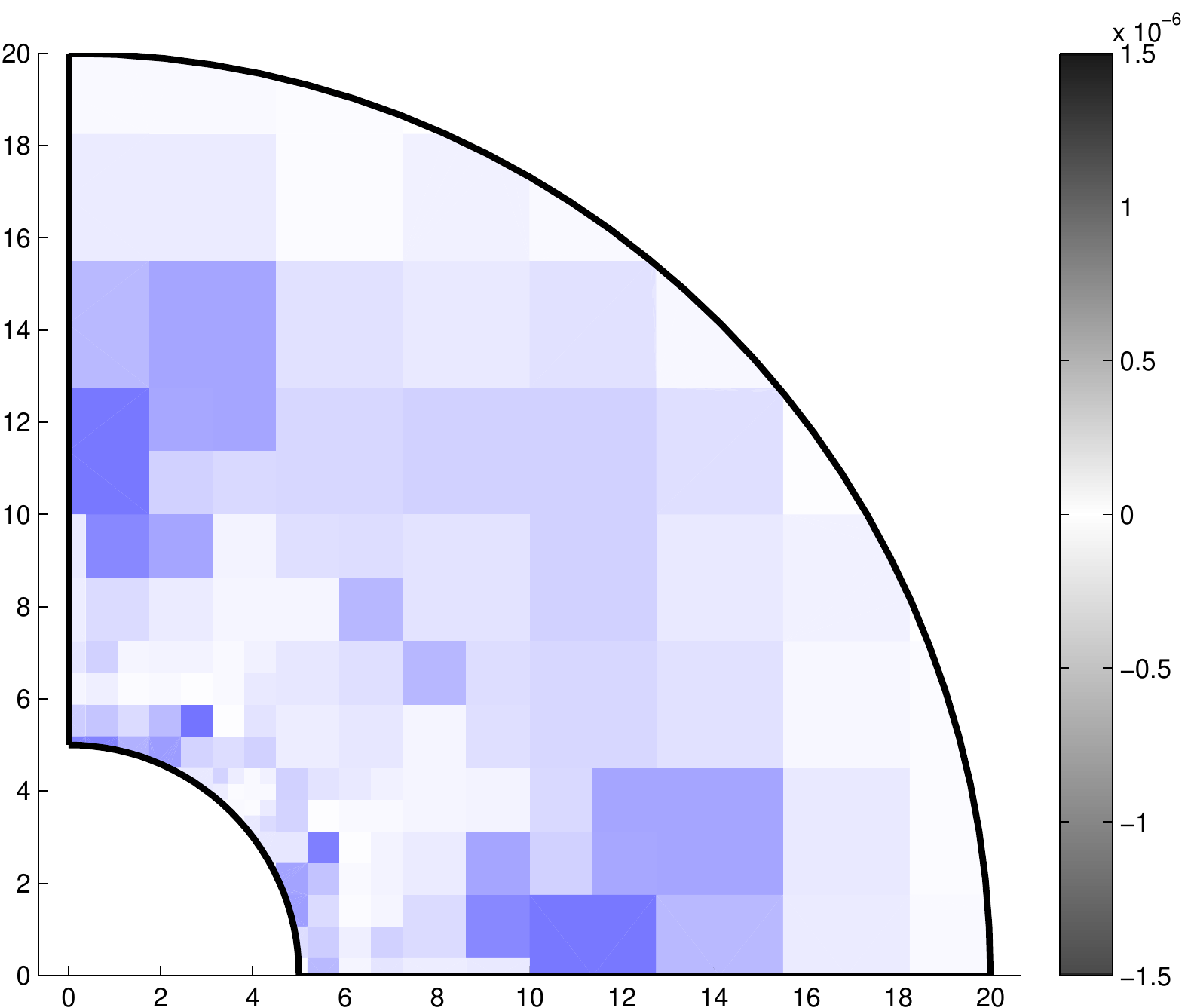} \\
   3 & \includegraphics[trim = 3.2cm 0.4cm 2.6cm 0.3cm, clip, width=0.4\textwidth]{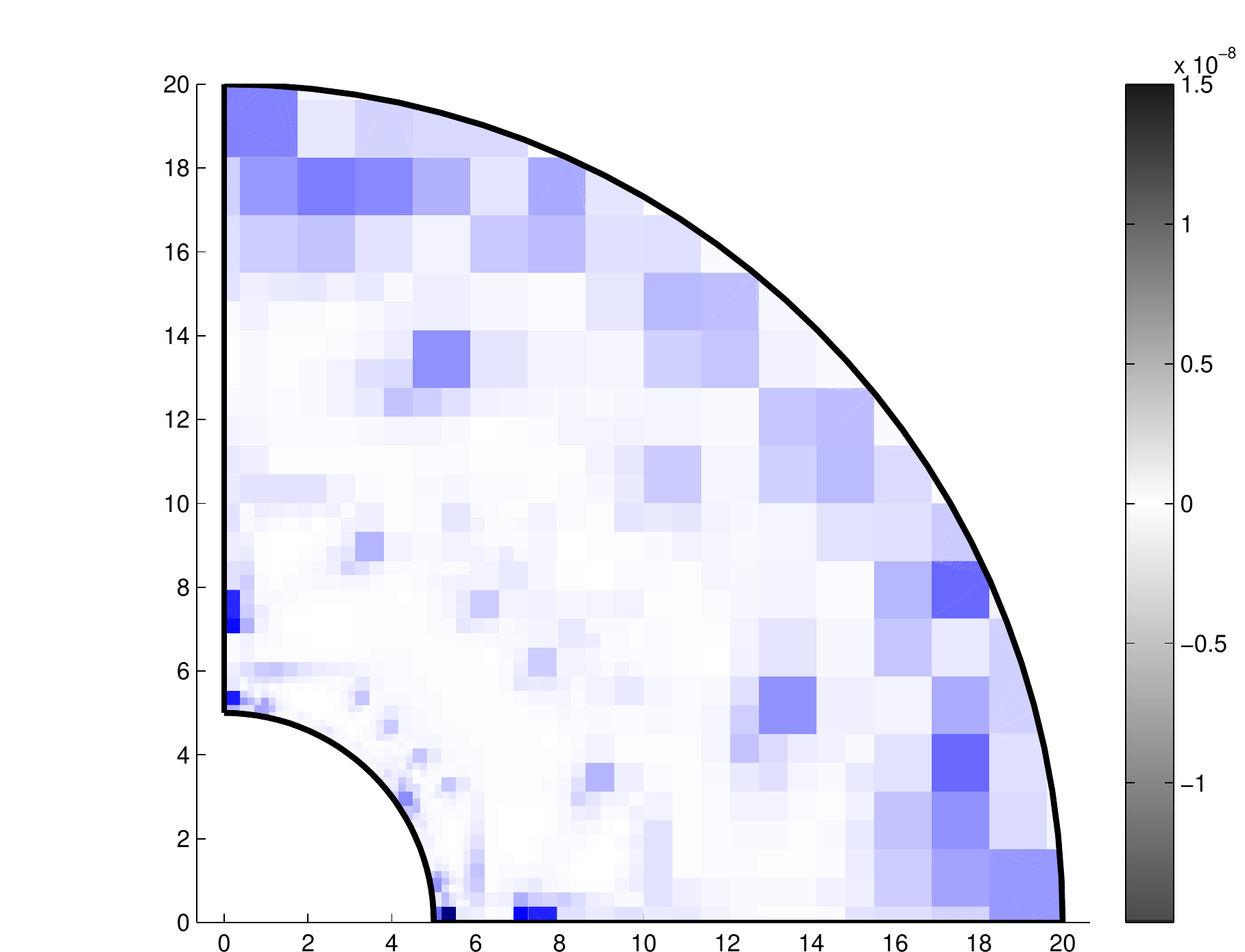} &   	\includegraphics[trim = 3.2cm 0.4cm 2.6cm 0.3cm, clip, width=0.4\textwidth]{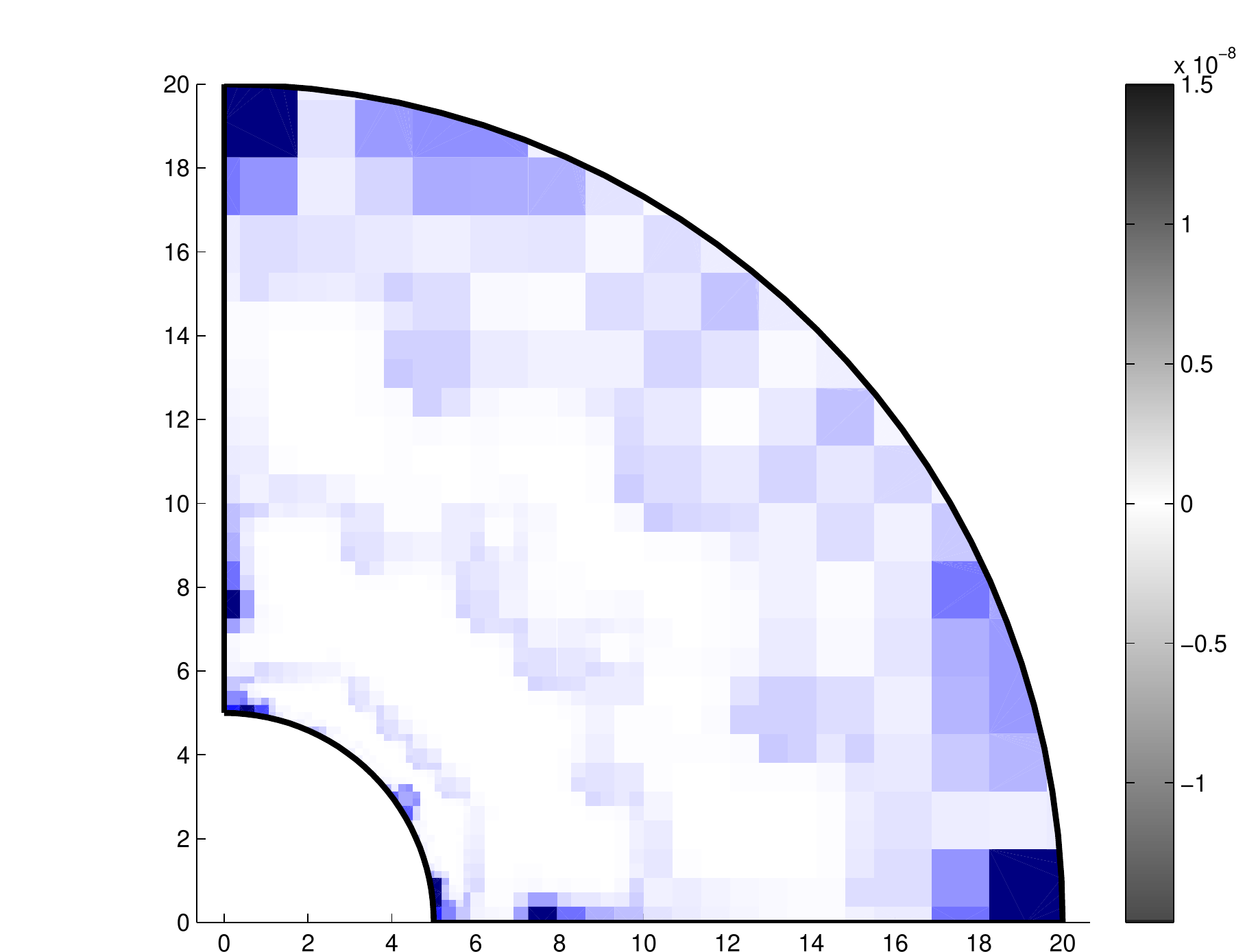} \\
   4 & \includegraphics[trim = 3.2cm 0.4cm 2.6cm 0.3cm, clip, width=0.4\textwidth]{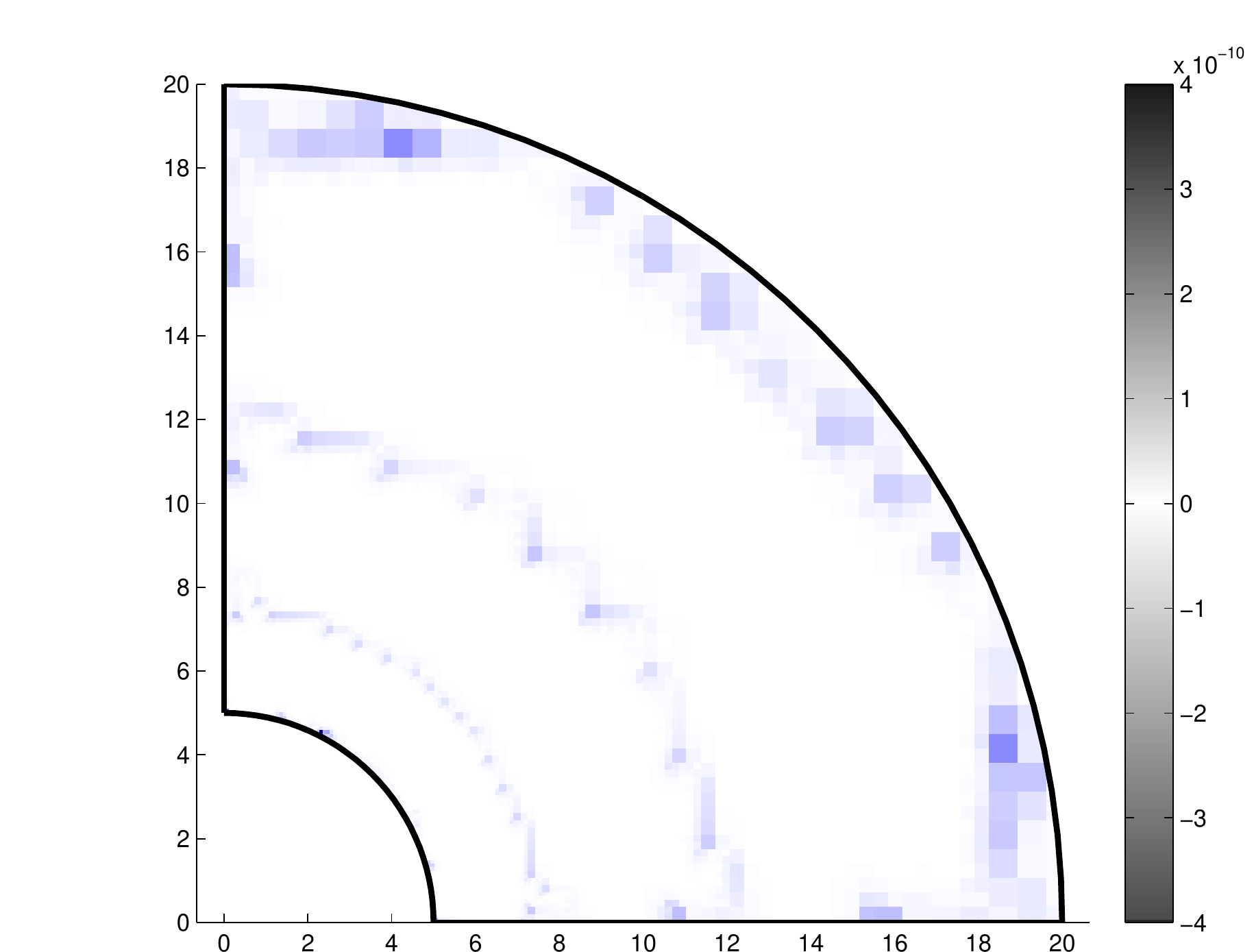} &   	\includegraphics[trim = 3.2cm 0.4cm 2.6cm 0.3cm, clip, width=0.4\textwidth]{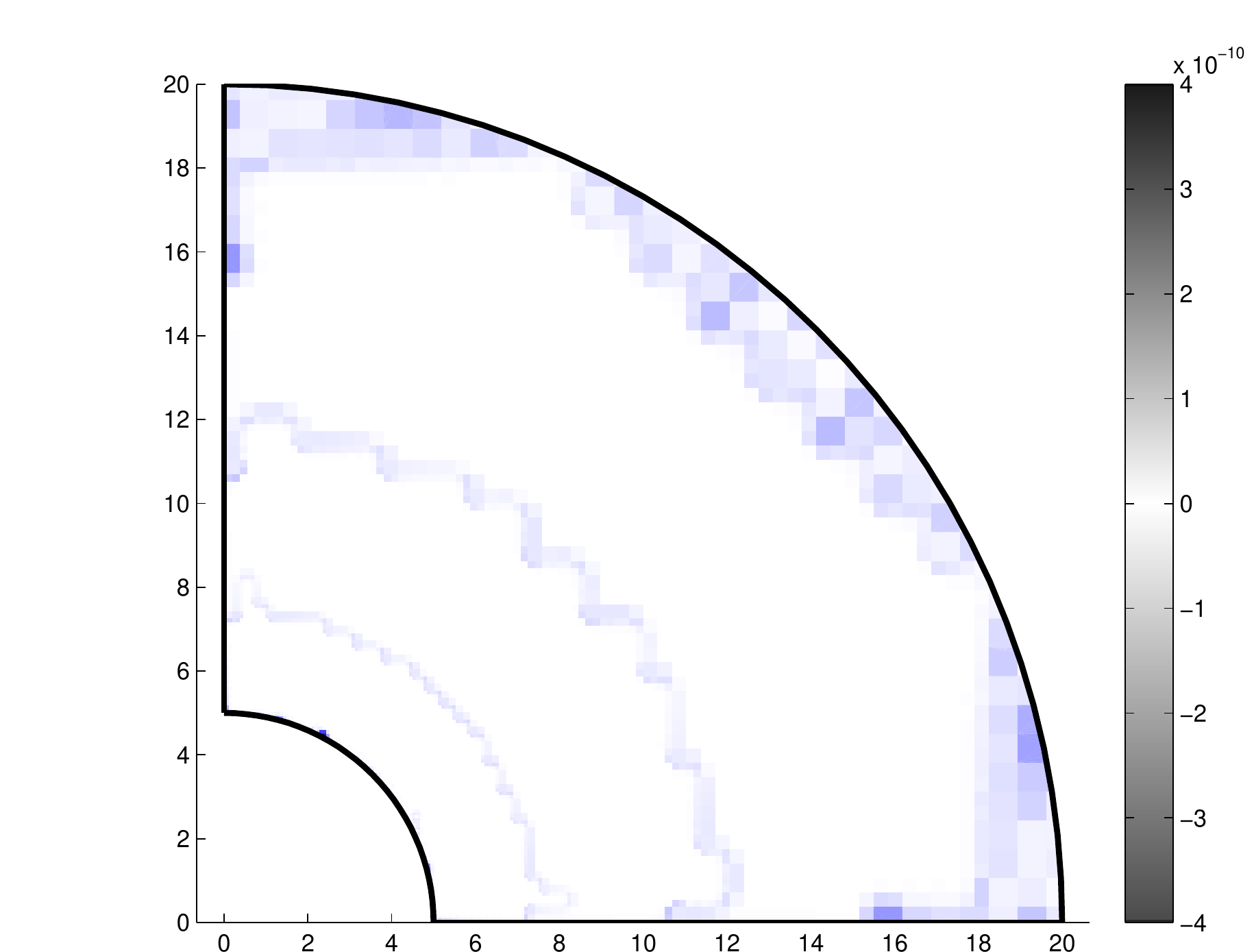} \\
   \end{tabular}
  \caption{Example 3. Pipe. Q4. Local exact error of the recovered solution (left),  local error estimates using $E^{*}_{3}$ estimator (right).}
  \label{fig:LocErrEst}
\end{figure}

For this problem we have also evaluated the distribution of the exact energy norm of the recovered solution $\left\|\be^*\right\|^2_k$ and the distribution of the error estimate $E^{*}_{3}|_k$ (subindex $k$ indicates that these errors are evaluated locally, at elements). In Figure \ref{fig:LocErrEst} we have represented these results for a sequence of \textit{h}-adapted meshes. We can observe that both error distributions are quite similar. These results show that the error estimator for the recovered solutions $E^{*}_{3}$ has a good behavior at global level, but also at local level.

%
%
\subsection{L-Shape. Singular problem}
\begin{figure}[!]
	\begin{tabular}{m{0.25\textwidth} l}

		\includegraphics[width=0.4\textwidth]{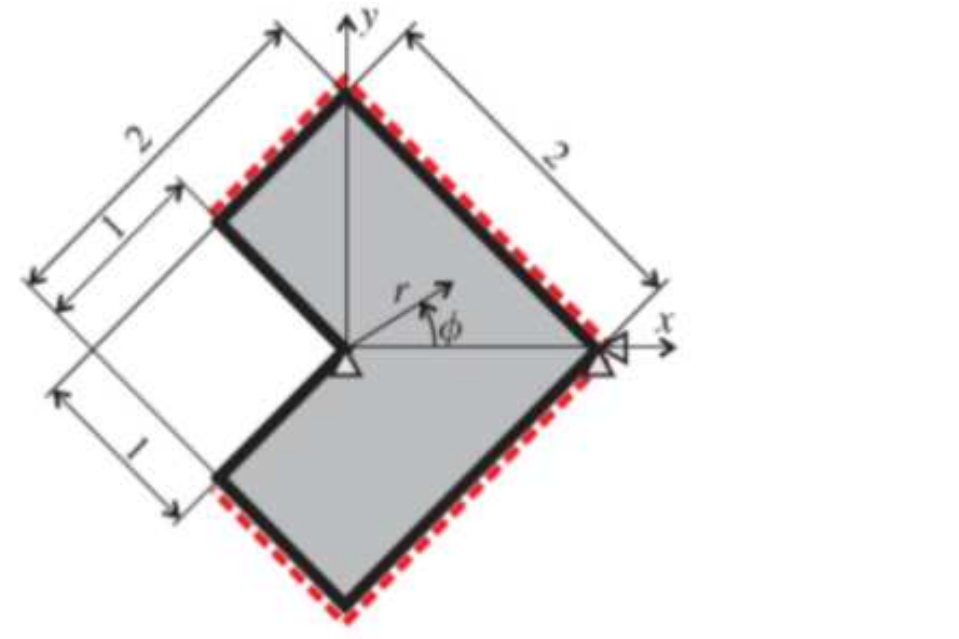} 

		&
		\parbox[l]{.75\textwidth}{
		\begin{align*}
			u_x(r,\upphi) &=\frac{1}{2G}r^L\left[\left(k-Q\left(L+1\right)\right)cos\left(L\upphi\right)-Lcos(\left(\left(L-2\right)\upphi\right)\right]\\
			u_y(r,\upphi) &=\frac{1}{2G}r^L\left[\left(k+Q\left(L+1\right)\right)sin\left(L\upphi\right)+Lcos(\left(\left(L-2\right)\upphi\right)\right]\\
			E&=1000 \; \nu=0.3 \\
			L&=0.544483736782464 \; Q=0.543075578836737  \\
			G&=\frac{E}{2(1+\nu)} \; k=3-4\nu
		\end{align*}
			}
	\end{tabular}
	\caption{Example 4: L-Shape model and analytical solution
		}
	\label{fig:LS_model}
\end{figure}

\begin{figure}[!]
\centering
\begin{tabular}{c c}
\begin{tikzpicture}
 	\begin{semilogxaxis}[
 	title={$\theta$},
 	ymax=1.5,
 	ymin=0,
 	xlabel={dof},
 	legend style={cells={anchor=west}, font=\small},
 	legend style={at={(0.98,0.02)},anchor=south east},
 	cycle list name=ageplot]
 	\draw (axis cs: 0,1)--(axis cs: 1e8,1); 
 	\addplot table[x=NDOF,y= Eff*_M1]{Data/LS-Q4_T012.csv};
 	\addplot table[x=NDOF,y= Eff*_M2]{Data/LS-Q4_T012.csv};
 	\addplot table[x=NDOF,y= Eff*_M3]{Data/LS-Q4_T012.csv};
 	\legend{{$E^*_1$},{$E^*_2$},{$E^*_3$}}
 	\end{semilogxaxis}
\end{tikzpicture}  
&
\begin{tikzpicture}
 	\begin{semilogxaxis}[
 	title={$\theta$},
 	ymax=1.5,
 	ymin=0,
 	xlabel={dof},
 	legend style={cells={anchor=west}, font=\small},
 	legend style={at={(0.98,0.02)},anchor=south east},
 	cycle list name=ageplot]
 	\draw (axis cs: 0,1)--(axis cs: 1e8,1);
 	\addplot table[x=NDOF,y= Eff*_M1]{Data/LS-Q8_T013.csv};
 	\addplot table[x=NDOF,y= Eff*_M2]{Data/LS-Q8_T013.csv};
 	\addplot table[x=NDOF,y= Eff*_M3]{Data/LS-Q8_T013.csv};
 	\legend{{$E^*_1$},{$E^*_2$},{$E^*_3$}}
 	\end{semilogxaxis}
 \end{tikzpicture} 
\end{tabular}
 \caption{Example 4: L-Shape. Effectivity index. Q4 (left) and Q8 (right).} 
\label{fig:LSeffect}
\end{figure}
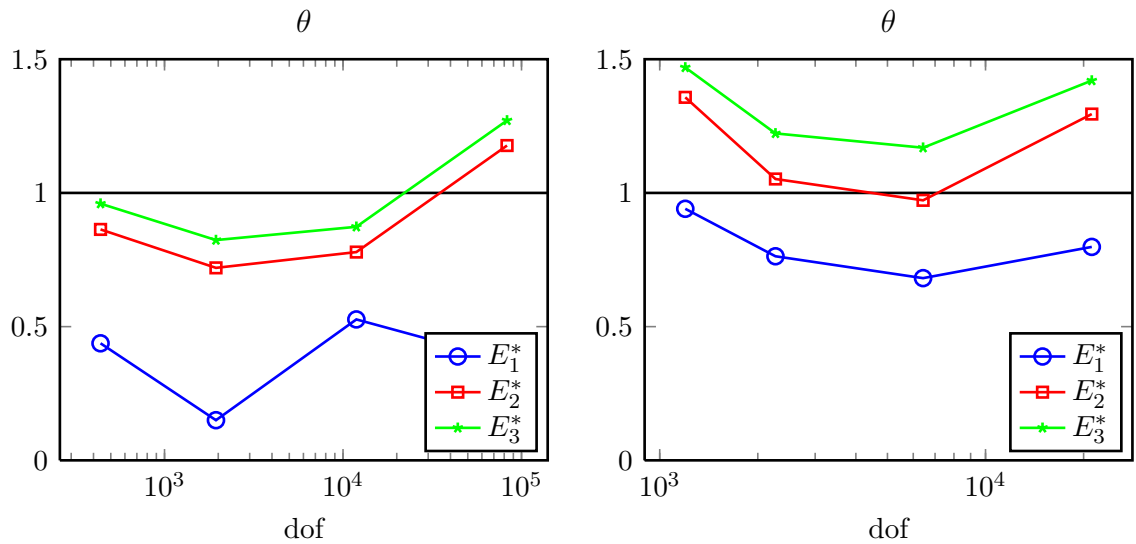

%
%
%
Consider the problem of an infinite plate with a V-notch subjected to tractions. We have considered the Mode I loading condition. The model of the problem, as defined in Figure \ref{fig:LS_model} has a singular solution in stresses. In this case we get similar results for the evolution of the global effectivity index  for estimators $E^{*}_{2}$ and $E^{*}_{3}$ (see Figure \ref{fig:LSeffect}), although $E^{*}_{3}$ exhibits the best results for linear elements. In the last meshes of the analysis we always obtain an increase of the effectivity index. This is because in the first 3 meshes of each analysis the mesh increases its density around the reentrant corner increasing the refinement level as we get closer to the singular point. The FE code can reach up to 20 refinement levels. In the last mesh the refinement level around the singularity would require higher refinement levels. The result is that, as higher refinement levels cannot be reached, we obtain an area around the singularity with elements of uniform size as opposite to graded meshes towards the singularity. This produces pollution errors and a decrease in the accuracy of the error estimation that leads to worse effectivity indexes.  

\begin{figure}[!]
  \centering
  \begin{tabular}{c c c}
   Mesh & $\left\|\be^*\right\|^2_k$ & $E^{*}_{3}|_k$ \\
   2 & \includegraphics[trim = 3.2cm 0.4cm 2.6cm 0.3cm, clip, width=0.4\textwidth]{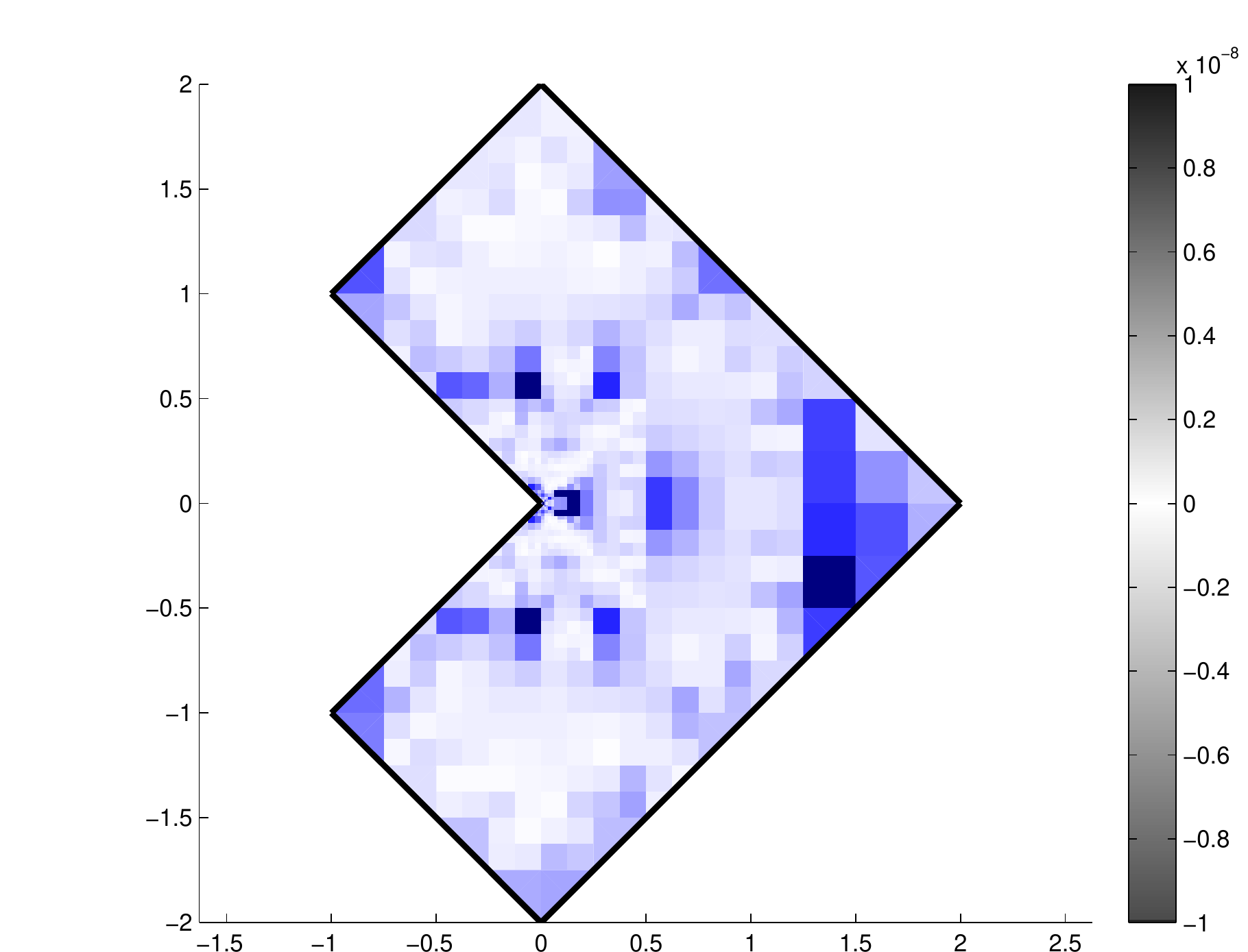} &   	\includegraphics[trim = 3.2cm 0.4cm 2.6cm 0.3cm, clip, width=0.4\textwidth]{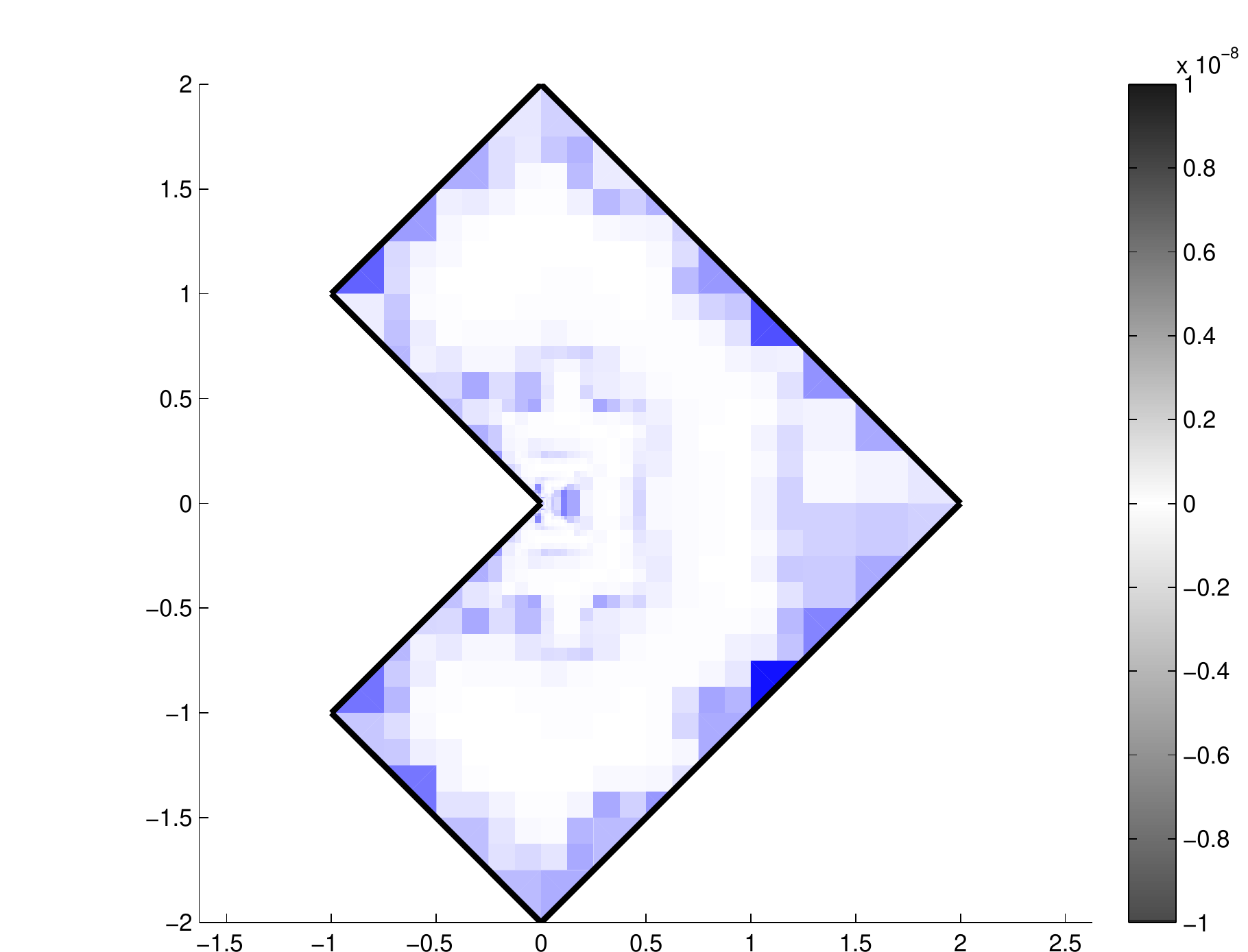} \\
   3 & \includegraphics[trim = 3.2cm 0.4cm 2.6cm 0.3cm, clip, width=0.4\textwidth]{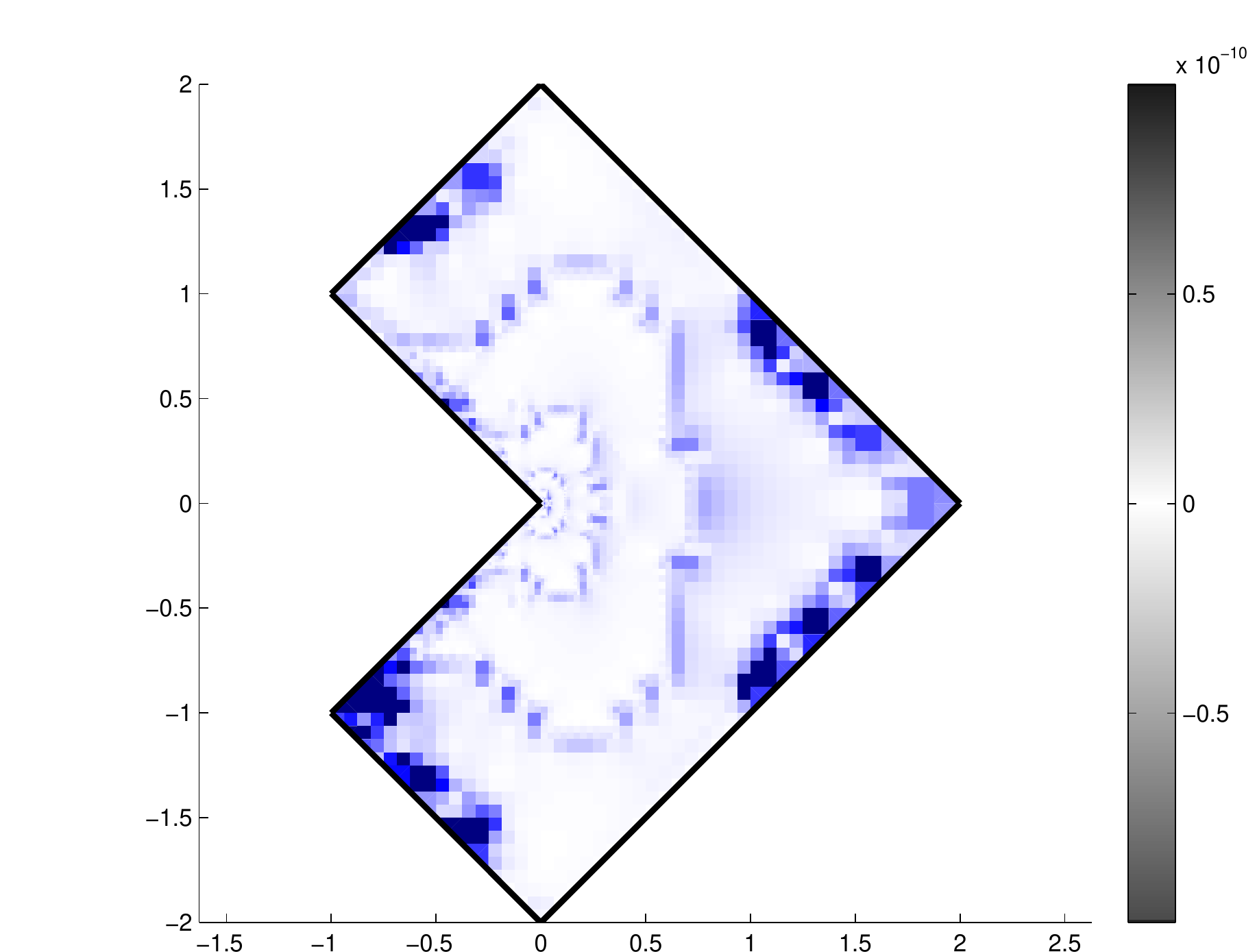} &   	\includegraphics[trim = 3.2cm 0.4cm 2.6cm 0.3cm, clip, width=0.4\textwidth]{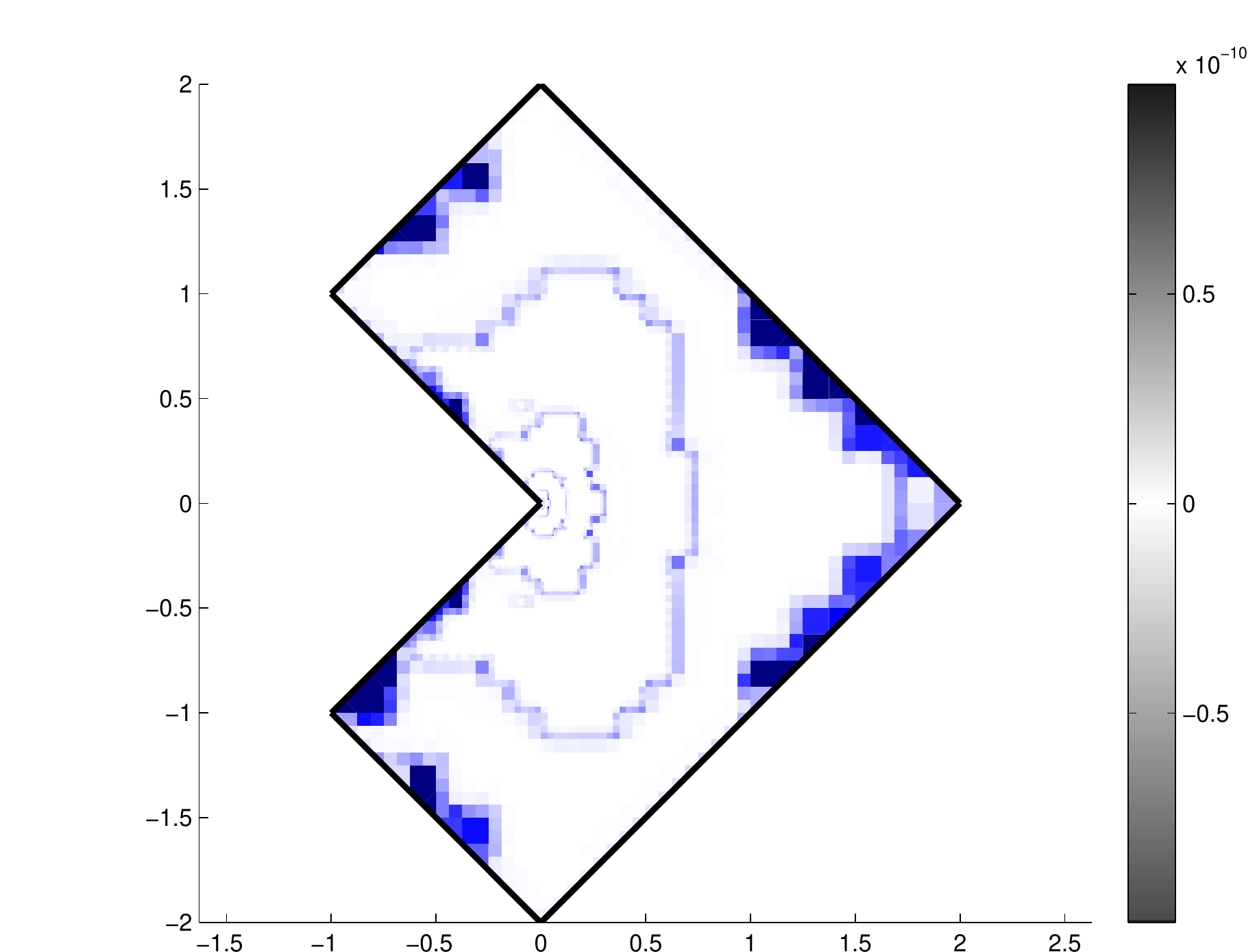} \\
   4 & \includegraphics[trim = 3.2cm 0.4cm 2.6cm 0.3cm, clip, width=0.4\textwidth]{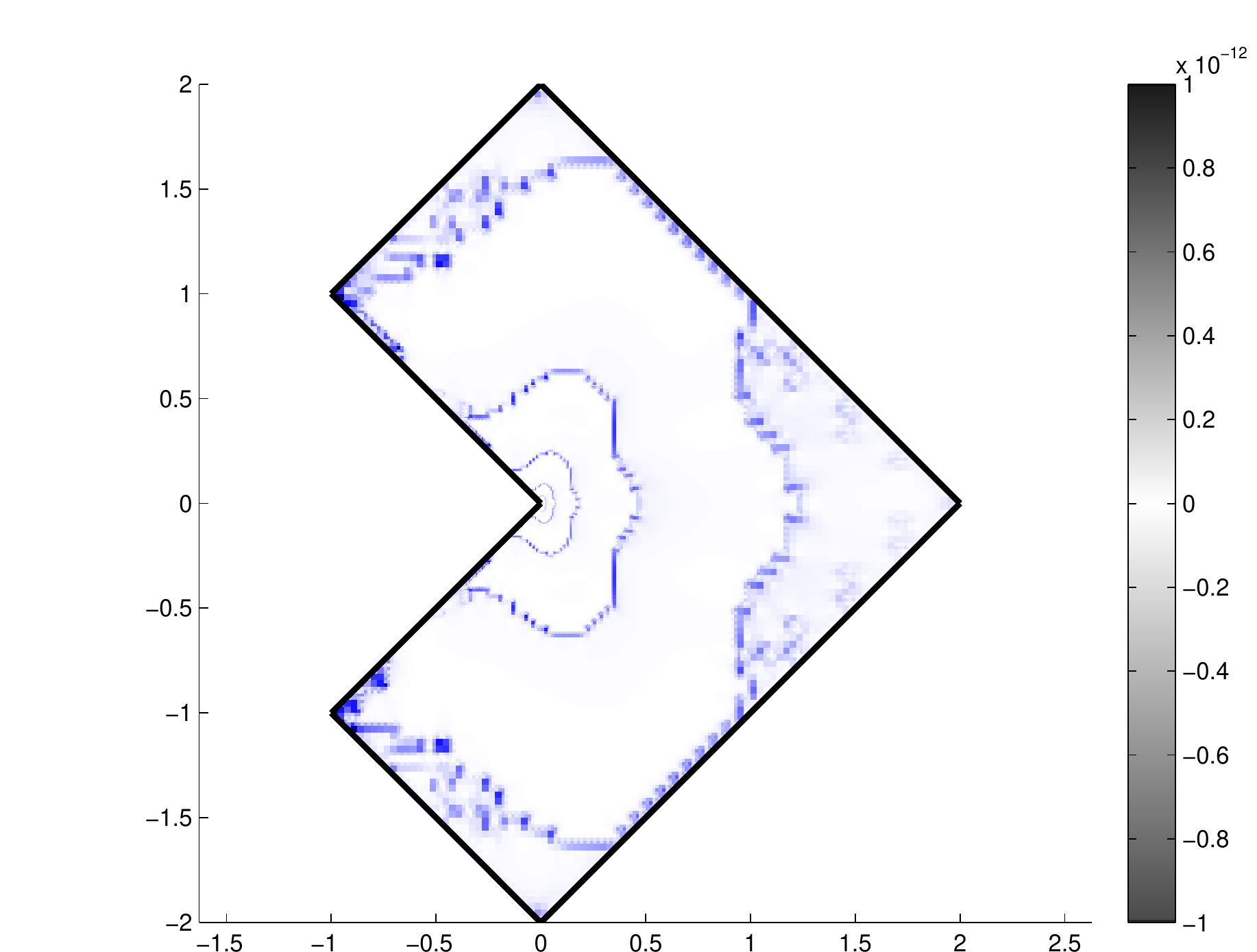} &   	\includegraphics[trim = 3.2cm 0.4cm 2.6cm 0.3cm, clip, width=0.4\textwidth]{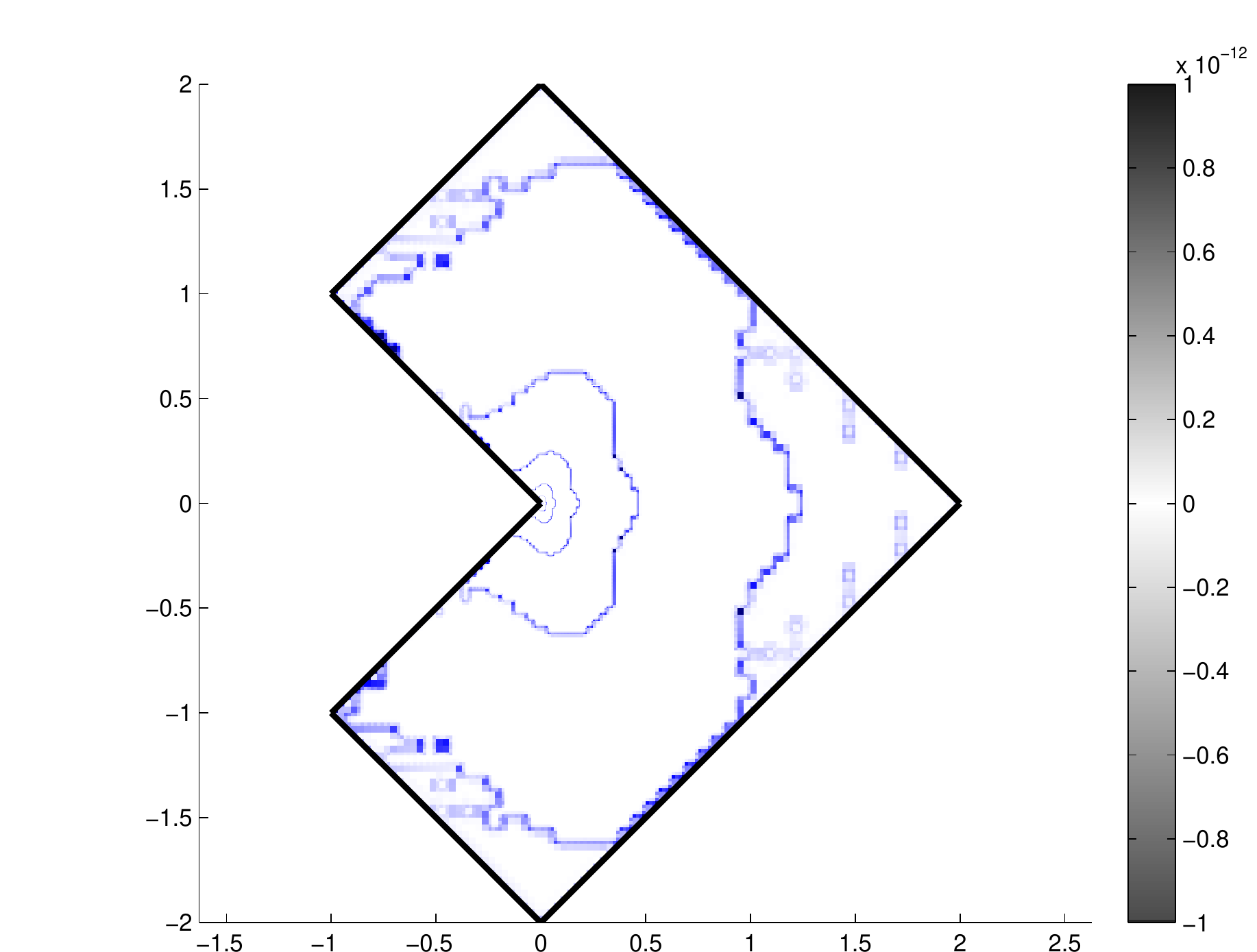} \\
   \end{tabular}
  \caption{Example 4. L-Shape. Q4. Local exact error of the recovered solution (left),  local error estimates using $E^{*}_{3}$ estimator (right).}
  \label{fig:LSLocErrEst}
\end{figure}

We show the local error index $E^{*}_{3}|_k$ and the exact error of the recovered solution, $\left\|\be^*\right\|^2_k$, in Figure \ref{fig:LSLocErrEst}. Observe that even for this singular problem the results are quite similar, concluding that, in this case $E^{*}_{3}|_k$ is a good error indicator.

%
%
\subsection{\textit{h}-adaptive process}
%
%
In previous sections we have shown that the method presented for the estimation of error of the recovered solution, using the estimator $E^*_3$, provides very accurate results with a very good global effectivity index. Now, we are going to show that this error estimation can be used in \textit{h}-adaptive processes.. 

During a \mbox{\textit{h}-adaptive} refinement process, using a standard FE compilation, the process reads as follows:

\begin{enumerate}
	\item Generate a FE mesh.
	\item Solve the FE problem.
	\item Estimate the error of the FE solution (locally and globally).
	\item If target error is smaller than the estimated error of the FE solution continue to step 5 else stop the process.
	\item Generate an \textit{h}-adapted FE mesh using the local FE error estimation
	\item Go to step 2.
\end{enumerate}

In this classical situation we are estimating the error of the raw FE solution, then we are using the raw FE solution $(\vm{u}^h,\vm{\sigma}^h)$ as output. However, when we use our recovery procedure we have already available an improved solution $(\vm{u}^*,\vm{\sigma}^*)$. So far, we were unable to estimate the error of this last solution, this output was not reliable. However, with the contribution presented in this report we have a way to obtain an accurate estimation of the error of this recovered solution: $E^*_3$, thus we can use $(\vm{u}^*,\vm{\sigma}^*)$ as the output of the analysis. The information about the error estimation in the recovered solution could be used in the \mbox{\textit{h}-adaptive} refinement process to obtain a solution $(\vm{u}^*,\vm{\sigma}^*)$ with the required accuracy:

\begin{enumerate}
	\item Generate a FE mesh.
	\item Solve the FE problem.
	\item Evaluate the local error estimate of FE solution $(\vm{u}^h,\vm{\sigma}^h)$.
	\item Evaluate the global value of the error estimate of the recovered solution $(\vm{u}^*,\vm{\sigma}^*)$ (costless procedure).
	\item If target error is smaller than the error of the recovered solution continue to step 5 else stop the process.
	\item Generate an \textit{h}-adapted FE mesh using the local FE error estimation 
	\item Go to step 2.
\end{enumerate}

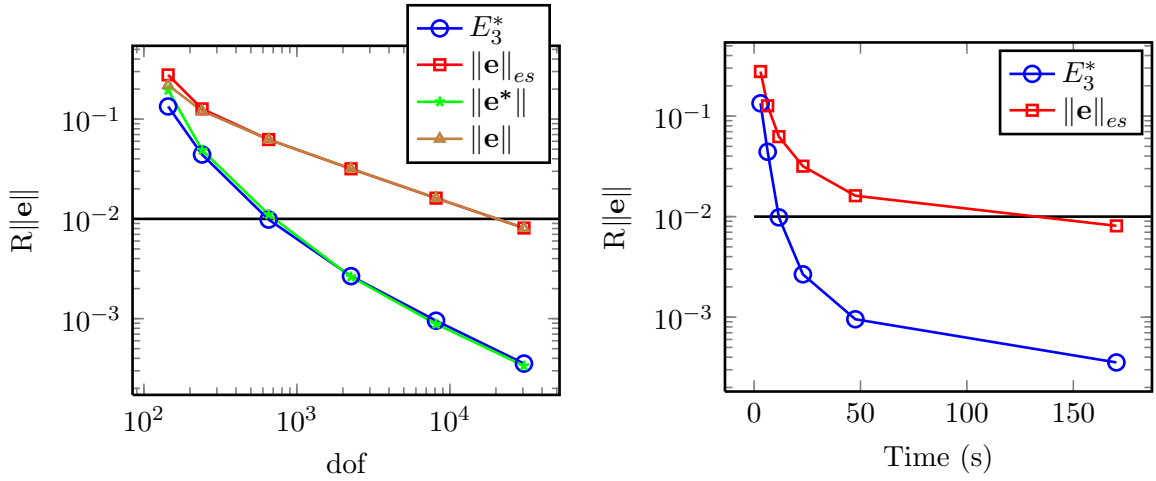
\begin{figure}[!]
\centering
	\begin{subfigure}[b]{0.45\textwidth}
	\centering
	\begin{tikzpicture}
	 	\begin{loglogaxis}[
	 	width=1\textwidth,
	 	xlabel={dof},
	 	ylabel={R$\left\|\vm{e}\right\|$},
	 	legend style={cells={anchor=west}, font=\small},
	 	legend style={at={(0.98,1.12)},anchor=north east},
	 	cycle list name=ageplot]
	 	\draw (axis cs: 0,0.01)--(axis cs: 1e8,0.01); 
	 	\addplot table[x=NDOF,y= R||e*es||]{Data/Cyl-Q4-T052.csv};
	 	\addplot table[x=NDOF,y= R||e_es||]{Data/Cyl-Q4-T052.csv};
	 	\addplot table[x=NDOF,y= R||e*||]{Data/Cyl-Q4-T052.csv};
	 	\addplot table[x=NDOF,y= R||e||]{Data/Cyl-Q4-T052.csv};
	 	\legend{{$E^*_3$},$\left\|\vm{e}\right\|_{es}$,   $\left\|\vm{e^*}\right\|$,$\left\|\vm{e}\right\|$}
	 	\end{loglogaxis}
	\end{tikzpicture} 
	\caption{Convergence analysis. Relative error in energy norm versus degree of freedom}
	\label{fig:ErrRec}
	\end{subfigure} \hspace{0.02\textwidth}
	\begin{subfigure}[b]{0.45\textwidth}
	\centering
	\begin{tikzpicture}
	 	\begin{semilogyaxis}[
	 	width=1\textwidth,
	 	xlabel={Time (s)},
	 	ylabel={R$\left\|\vm{e}\right\|$},
	 	legend style={cells={anchor=west}, font=\small},
	 	legend style={at={(0.98,0.98)},anchor=north east},
	 	cycle list name=ageplot]
	 	\draw (axis cs: 0,0.01)--(axis cs: 1e3,0.01); 
	 	\addplot table[x=Acc Time,y=R||e*es||]{Data/Cyl-Q4-T052.csv};
	 	\addplot table[x=Acc Time,y=R||e_es||]{Data/Cyl-Q4-T052.csv};
	 	\legend{{$E^*_3$},$\left\|\vm{e}\right\|_{es}$}
	 	\end{semilogyaxis}
	\end{tikzpicture} 
	\caption{Computational cost analysis. Relative error in energy norm versus time.}
	\label{fig:TimeRec}
	\end{subfigure}
\caption{Example 3. \textit{h}-adaptive analysis with Q4 elements. The black line represents the prescribed relative error in energy norm (1\%).}
\label{fig:RecSoluc}
\end{figure}

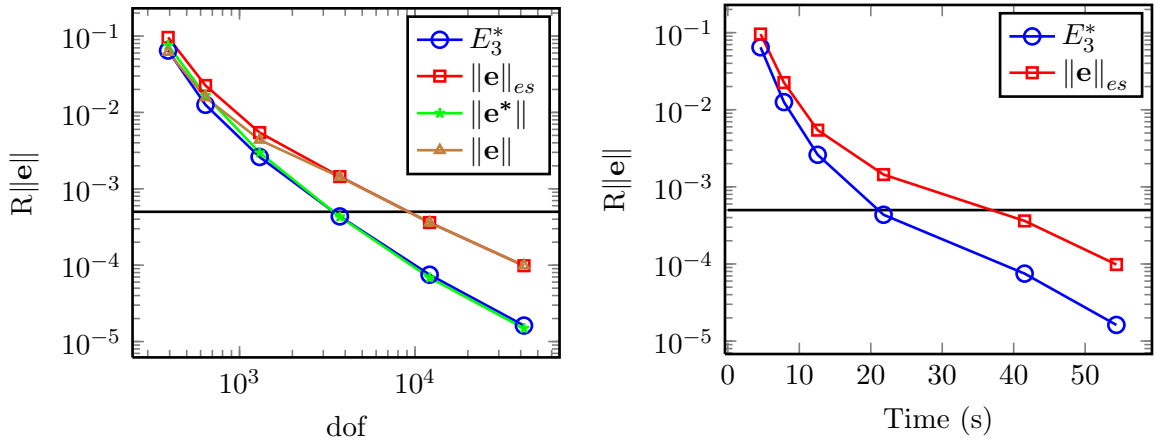
\begin{figure}[!]
\centering
	\begin{subfigure}[b]{0.45\textwidth}
	\centering
	\begin{tikzpicture}
	 	\begin{loglogaxis}[
	 	width=1\textwidth,
	 	xlabel={dof},
	 	ylabel={R$\left\|\vm{e}\right\|$},
	 	legend style={cells={anchor=west}, font=\small},
	 	legend style={at={(0.98,0.98)},anchor=north east},
	 	cycle list name=ageplot]
	 	\draw (axis cs: 0,0.0005)--(axis cs: 1e8,0.0005); 
	 	\addplot table[x=NDOF,y= R||e*es||]{Data/Cyl-Q8-T053.csv};
	 	\addplot table[x=NDOF,y= R||e_es||]{Data/Cyl-Q8-T053.csv};
	 	\addplot table[x=NDOF,y= R||e*||]{Data/Cyl-Q8-T053.csv};
	 	\addplot table[x=NDOF,y= R||e||]{Data/Cyl-Q8-T053.csv};
	 	\legend{{$E^*_3$},$\left\|\vm{e}\right\|_{es}$,   $\left\|\vm{e^*}\right\|$,$\left\|\vm{e}\right\|$}
	 	\end{loglogaxis}
	\end{tikzpicture} 
	\caption{Convergence analysis. Relative error in energy norm versus degree of freedom}
	\label{fig:ErrRec8}
	\end{subfigure} \hspace{0.02\textwidth}
	\begin{subfigure}[b]{0.45\textwidth}
	\centering
	\begin{tikzpicture}
	 	\begin{semilogyaxis}[
	 	width=1\textwidth,
	 	xlabel={Time (s)},
	 	ylabel={R$\left\|\vm{e}\right\|$},
	 	legend style={cells={anchor=west}, font=\small},
	 	legend style={at={(0.98,0.98)},anchor=north east},
	 	cycle list name=ageplot]
	 	\draw (axis cs: 0,0.0005)--(axis cs: 1e3,0.0005); 
	 	\addplot table[x=Acc Time,y=R||e*es||]{Data/Cyl-Q8-T053.csv};
	 	\addplot table[x=Acc Time,y=R||e_es||]{Data/Cyl-Q8-T053.csv};
	 	\legend{{$E^*_3$},$\left\|\vm{e}\right\|_{es}$}
	 	\end{semilogyaxis}
	\end{tikzpicture} 
	\caption{Computational cost analysis. Relative error in energy norm versus time.}
	\label{fig:TimeRec8}
	\end{subfigure}
\caption{Example 3. \textit{h}-adaptive analysis with Q8 elements. The black line represents the prescribed relative error in energy norm (0.05\%).}
\label{fig:RecSoluc8}
\end{figure}

The \textit{h}-adaptive procedure is guided in both cases by the well established techniques based on the error estimation of the FE solution, which will simultaneously decrease the error of the FE solution and the error of the recovered solution. Note that the main difference between both processes is, simply, the stopping criterion. In this new approach the process stops when the estimated error of the recovered solution is smalled than the target error. The error in energy norm of the recovered solution is, in general, lower than the error in the FE solution and the convergence rate of the error in energy norm of the recovered solution is higher than that of the FE solution, thus, the recovered solution reaches the prescribed error level with less degrees of freedom than the FE solution. Therefore, this second method produces important savings in total computational cost of the analysis. Figures \ref{fig:RecSoluc} and \ref{fig:RecSoluc8} show the evolution of the errors during the \mbox{\textit{h}-adaptive} refinement process and the computational cost needed to obtain a certain accuracy level, for Example 3 (cylinder problem) for Q4 and Q8 elements.

The horizontal black lines in figures \ref{fig:RecSoluc} and \ref{fig:RecSoluc8} represent the error level of the solution prescribed by the analyst. Red and brown lines represent the error (exact and estimated) of the standard FE output $(\vm{u}^h,\vm{\sigma}^h)$. Blue and green lines represent the error (exact and estimated) of the recovered field $(\vm{u}^*,\vm{\sigma}^*)$. We can observe that the error estimation of recovered solution accurately represents the exact error for both linear (Q4) and quadratic (Q8) elements. 


The figures show a considerable improvement of the \textit{h}-adaptive process to reach a prescribed error level. Note that, for Q4 elements, to reach the prescribed relative error in energy norm (1\%), the standard \textit{h}-adaptive strategy based on the accuracy of the FE solution would stop the process after an analysis with 30492 dofs and 170.2 s whereas with the proposed \textit{h}-adaptive method based on the accuracy of the recovered solution the process would stop with only 654 dofs and 11.7 s. With Q8 elements reaching a 0.05\% prescribed error, the standard \textit{h}-adaptive procedure would stop the process after an analysis wiht 41724 dofs and 54.4 s, whereas the proposed method would only require 3728 dofs and 21.8 s.

%
%
\section{Conclusions}
%
%
Error estimation of the recovered solution could be very useful to decrease the computational cost of FE analysis. Error estimator $E^*_3$ efficiently predicts the error in energy norm of the recovered solution both at local and global levels. This accuracy is not mathematically justified. The results obtained definitively require the support of mathematical proofs.

\section{Acknowledgements}

This work was supported by the EPSRC grant EP/G042705/1 ``Increased Reliability for Industrially Relevant Automatic Crack Growth Simulation with the eXtended Finite Element Method''. 

This work has received partial support from the research project DPI2010-20542 of the Ministerio de Econom\'ia y Competitividad (Spain). The financial support of the FPU program (AP2008-01086), the funding from Universitat Polit\`{e}cnica de Val\`{e}ncia and Generalitat Valenciana (PROMETEO/2012/023) are also acknowledged. 

All authors also thank the partial support of the Framework Programme 7 Initial Training Network Funding under grant number 289361 ``Integrating Numerical Simulation and Geometric Design Technology."

\clearpage
\bibliographystyle{wileyj-nourldoi}
\bibliography{library}

\end{document}